\documentclass[aap,MSNbibl,nameyear,nochecklpage,dvips]{arximspdf}
\usepackage{mathbh}

\doi{10.1214/09-AAP640}
\volume{20}
\issue{3}
\pubyear{2010}
\firstpage{907}
\lastpage{950}

\makeatletter
\newtheorem{theorem}{Theorem}[section] %
\newtheorem{cor}[theorem]{Corollary}
\newproclaim{Example}[theorem]{Example}
\newtheorem{lemma}[theorem]{Lemma}
\newtheorem{propos}[theorem]{Proposition}
\newproclaim{remark}[theorem]{Remark}
\newproclaim{Remark}[theorem]{Remark}

\newcommand{\refeq}[1]{(\ref{#1})}
\makeatother

\begin{document}
\begin{frontmatter}

\title{A functional limit theorem for\\ the profile of $b$-ary trees}
\runtitle{Functional limit theorem for the profile of $b$-ary trees}

\begin{aug}
\author[A]{\fnms{Eva-Maria} \snm{Schopp}\ead[label=e1]{eva.schopp@gmx.de}\corref{}}
\runauthor{E.-M. Schopp}
\affiliation{University of Freiburg}
\address[A]{Department of Mathematical Stochastics\\
University of Freiburg\\
Eckerstr. 1\\
79104 Freiburg\\
Germany\\
\printead{e1}}
\end{aug}

\received{\smonth{7} \syear{2008}}
\revised{\smonth{8} \syear{2009}}

%
\begin{abstract}
In this paper we prove a functional limit theorem for the weighted
profile of a $b$-ary tree.
For the proof we use classical martingales connected to branching
Markov processes and a generalized version of the profile-polynomial martingale.
By embedding, choosing weights and a branch factor in a right way, we
finally rediscover the profiles of some well-known discrete time trees.
\end{abstract}

%


\end{frontmatter}
%

\section{Introduction}\label{intr}
The profile is the set or sequence of numbers of nodes at each level of
a tree. It is a fine tree shaped parameter related to many other
important shape characteristics such as the total path length (the sum
of the distances of all nodes to the root), depth (the distance of a
random node to the root), height (the maximal distance of a node to the
root), saturation level (the minimal distance of an external node to
the root) and width (number of nodes at the most abundant level).

In general, we distinguish between two main groups of trees. In the
first group we collect all the trees by their height which increases
through the square root of the number of nodes. The prime examples of
these kind of trees are Galton--Watson trees conditioned on the total
progeny or simply generated trees.
It is shown in \citet{Aldous93} that the simply generated trees
studied by \citet{MeirMoon78} and conditioned Galton--Watson trees are
the same.
Moreover, the profile of conditioned Galton--Watson trees is further
investigated in \citet{DrmotaGittenberger97}, \citet{Kersting98} and
\citet{Pitman99b} with further references.

In the second group we collect all the trees with logarithmically
growing height. In contrast to conditioned Galton--Watson trees, the
profiles of these trees have received less attention in the past.
However, there has been much done on this topic in the last few years.
The methods used to derive limit theorems for normalized profiles range
from the method of moments, the contraction method to analytical tools,
including saddlepoint methods, Mellin transforms, Poissonization,
de-Poissonization, singularity analysis, application of generating
functions and uniform asymptotic analysis. We refer to
\citet{AldousShields88}, Drmota and Hwang (\citeyear{DrmotaHwang05b,DrmotaHwang05a}), \citet
{FuchsHwangNeininger05}, \citet{Hwang05}, \citet{DevroyeHwang06},
\citet{DrmotaJansonNeininger06}, \citet{Park06} and \citet{ParkHwangNicodemeSzpankowski06}.

Furthermore, \citet{ChauvinDrmotaJabbour01} and Chauvin et al.~(\citeyear
{ChauvinKleinMarckertRouault05}) used martingale methods to obtain
limit theorems for the profile of the random binary search tree.
It is a classical result of \citet{Jabbour-Hattab01} that the
profile-polynom of the binary search tree,\vspace*{1pt} $ M_n(z):=\frac
{1}{C_n(z)}\sum_{l}U_l(n)z^l$, is a discrete
time martingale where $U_l(n)$ denotes the number of nodes in
generation $l$ of a random binary search tree of size $n$ [cf. \refeq
{eqU_l(n)}] and $C_n(z)$ is a normalizing deterministic constant for
$z \in\mathbf{\mathbb{C}}$ fixed [see \refeq{eqC_n} for a definition]. With
this in mind, classical martingale convergence results may now be
applied to investigate the (asymptotic) behavior of this martingale
and the asymptotics of the profile of the binary search tree [see
\citet{ChauvinDrmotaJabbour01}]. Instead of analyzing $M_n(z)$ directly,
\citet{ChauvinKleinMarckertRouault05} showed that the discrete time
martingale $M_n(z)$ is deeply related to
the well-studied classical Yule-time martingale $M(t,z)$ of the
corresponding continuous time tree, the Yule tree [cf. \refeq{M(t,z)}].
Finally, they strengthened classical convergence results in order to
relate the uniform convergence of $M(t,z)$ in compact sets to
corresponding uniform asymptotic results of $M_n(z)$.

Apart from the profile of the binary search tree, the profile of the
random recursive tree has been studied with martingale methods [cf.
\citet{DrmotaHwang05b}]. Drmota and Hwang (\citeyear{DrmotaHwang05a})
showed that the profile is not concentrated around the mean. Their
proof is based on showing that the variance undergoes four phase
transitions and exhibits a bimodal behavior in contrast to the
unimodality of the expected value of the profiles. As a consequence,
the profile is not concentrated around the mean. For example, around
the most numerous level (where the width is attained) the variance is
small [$O(\log n)$] and the profile is concentrated with a Gaussian
limit, but $\log n$ away from this level the variance increases, and
there is no concentration anymore
[cf. Drmota and Hwang (\citeyear{DrmotaHwang05b,DrmotaHwang05a}), and further remarks in Example \ref{ex3.3.2}].
Recently, \citet{Sulzbach08} used the martingale method of
\citet{ChauvinKleinMarckertRouault05} to show a limit
theorem for the profile of plane oriented recursive trees.

Obviously not all trees have a martingale structure,
and, furthermore, the method rapidly becomes costly when we leave the
well-known binary search tree and recursive tree cases. The last
statement is based on the fact that corresponding statements about
other trees, and especially their continuous time matches, are not
available which differs from assertions about the well-studied Yule
tree process.

However, there is a large class of trees for which the martingale
method gives uniform convergence results for the profile in any
suitable compact set. In this paper we are interested in the kind of
trees that can be studied with martingale methods. We want to find out
how the method of \citet{ChauvinKleinMarckertRouault05} can be
generalized to derive asymptotic results for other trees besides the
ones that have been studied so far.

\citet{BroutinDevroye05} constructed a class of
continuous time edge-weighted trees and proved, using large deviation
techniques, a
general law of large numbers for the height of these trees.
This class of edge-weighted trees leads, when the tree process is
stopped correctly, to various well-known discrete time trees as, for
example, to random binary search trees and to random recursive trees.
The initial point of their model is the construction of an infinite,
$b$-ary tree where each node $u$ in that tree is assigned, independently,
a random vector,
\[
\bigl(\bigl(Z_1^{(u)},E_1^{(u)}\bigr), \bigl(Z_2^{(u)}, E_2^{(u)}\bigr), \ldots,
\bigl(Z_b^{(u)},E_b^{(u)}\bigr)\bigr),
\]
where, for instance, $E_j^{(u)}$ is the lifetime of $u$'s $j$th child,
and $Z_j^{(u)}$ is a specific weight assigned to the $j$th child of $u$
[cf. \refeq{vector}].
Let Exp$(\mu)$ denote the exponential distribution with the first
moment equal to $\frac{1}{\mu}, \mu>0$, in other words, the
distribution with Lebesgue density $f(x)=\mu e^{-\mu x}\mathbh{1}_{[0,\infty
)}(x)$ for $ x \in\mathbb{R}$.
A continuous time tree process is Markovian with a classical martingale
structure if and only if the lifetimes are independent of each other
and exponentially distributed (see, e.g., Harris [(\citeyear{Harris63}), Chapter V.2],
\citet{Watanabe67}, \citet{JoffeLeCamNeveu73}, Athreya and Ney [(\citeyear{AthreyaNey75}), Chapter
III], \citet{Kingman75}, \citet{Biggins77b}, \citet{Wang80}, \citet{Uchiyama82}, \citet{Neveu87}, \citet{Biggins91}
and \citet{Biggins92}). Essentially, this follows from the memoryless
property of the exponential distribution.

Our purpose is now to show that these trees have the right martingale
properties in order to generalize the method of \citet{ChauvinKleinMarckertRouault05}.
The main theorem of this paper, Theorem \ref{theo2.3.9}, states that
the normalized profile of a $b$-ary weighted tree converges almost
surely to a limit as the number of nodes in the tree tends to infinity.
This limit can be identified as the almost sure limit of a discrete
time martingale and, also, as the unique solution of expectation one of
a fixed point equation.

The plan for the rest of the paper is the following:
First, in Section \ref{section6.1} we will introduce the tree model of
\citet{BroutinDevroye05} and define a closely related branching Markov
process. Second, in Section \ref{section6.2} we will define the corresponding
continuous time martingale associated with the branching Markov process
of Section \ref{section6.1} and its discrete time analog in Section
\ref{section6.3}.
Next, in Section \ref{section6.4} we will elaborate the relationship
between those two martingales and formulate our main result, Theorem
\ref{theo2.3.9}
in Section \ref{section6.5}.
Finally, in Section \ref{section6.6} we will show the applicability of
our main Theorem \ref{theo2.3.9} based on some examples as, for
instance, the well-known random binary search tree, the random
recursive tree, random lopsided trees and random plane oriented trees.


\section{The framework}\label{section6.1}
In this section we describe the tree model of \citet
{BroutinDevroye05}. Let $\hat{T}_{\infty}$ be an infinite, complete
$b$-ary tree with $b\geq2$.
We assign to each node a label,
%
\begin{equation}\label{eqlabels}
u \in\mathcal{U}:=\{\o\}\cup\bigcup_{n=1}^{\infty}\{1, \ldots,
b\}^n,
\end{equation}
in the natural way; the root node, which will be denoted by \o, has $b$
children which are called $1, \ldots, b$. In the same manner, every
node $u$ has children denoted by $u1,\ldots, ub$. Generally, if
$u=u_1\cdots u_l$ is a node, and $v:=v_1\cdots v_k$ is a sequence with
$v_j \in\{1,\ldots, b\}, j=1,\ldots, k$, we set
\[
uv=u_1\cdots u_l v_1\cdots v_k
\]
and call $u$ an ancestor of $uv$.

For each node $u$ we create independently a random $b$-vector,
%
\begin{equation}\label{vector}
\bigl(\bigl(Z_1^{(u)},E_1^{(u)}\bigr), \bigl(Z_2^{(u)}, E_2^{(u)}\bigr), \ldots,
\bigl(Z_b^{(u)},E_b^{(u)}\bigr)\bigr),
\end{equation}
where for $j=1,\ldots, b,$ $(Z_j^{(u)},E_j^{(u)})$\vspace*{-1pt} is the vector
assigned to the edge from node~$u$ to its $j$th child.\vspace*{-1pt}
Here, $Z_j^{(u)}$ represents the weight of the edge from node $u$ to
its $j$th child, and $E_j^{(u)}$ is the lifetime of node $u$'s $j$th child.

Each couple $(Z_j^{(u)}, E_j^{(u)})$ is distributed as $(Z,E)$ for
independent $Z$ and $E$.
Note that the lifetime $E_j^{(u)}$ is independent of $E_r^{(v)}$ for
different $u$ and $v$ or different $j$ and~$r$.
All $E$'s are independent of any $Z$'s, but we allow dependence of
$Z_1^{(u)}, \ldots, Z_b^{(u)}$ for any node $u$ in the tree $\hat
{T}_{\infty}$.

We also assume that $E$ is exponentially distributed with mean one to
ensure Markov properties and that $Z$ is a lattice distribution with
values in $\mathbb{Z}^d$ for some $d\in\mathbb{N}$.

\begin{Remark}
In \citet{BroutinDevroye05} it is assumed that $E$ and $Z$ are
nonnegative independent random variables where $E$ is not mono-atomic,
has no atom at zero and the following property holds:
\[
\inf\{x\dvtx P(E>x)>0\}=0.
\]
As mentioned in their concluding remark, their model and their proof of
the law of large numbers for the height can be extended to more general
cases allowing dependence of $E$ and $Z$ as well as multi-dimensional versions.
In \citet{BroutinDevroyeMcLeish07} the assumption that the components
of the random vectors attained to each node are independent is skipped
in order to obtain
further height results of, for example, pebbled trees and others.

Let $\pi(u)$ be the set of edges from the root to node $u$. Further,
let $(Z_e, E_e)$ be the couple of random variables assigned to edge
$e$, and let
%
\begin{equation}\label{hat_T_t}
\hat{T}_t:=\biggl\{u \in\hat{T}_{\infty}\dvtx  G_u:=\sum_{e \in\pi(u)}E_e
\leq t\biggr\}
\end{equation}
be the subtree of $\hat{T}_{\infty}$ consisting of the nodes that
deceased before time $t$.
We are interested in the external profile of $\hat{T}_t$,
%
\begin{equation}\label{hat_rho}
\hat{\rho}_t(l):= \biggl| \biggl\{ u \in\partial\hat{T}_t:  D_u:=\sum
_{e \in\pi(u)}Z_e =l \biggr\} \biggr|
\end{equation}
for $t \in\mathbb{R}_+$ and $l \in\mathbb{Z}^d$,
where
\begin{eqnarray*}
&&\partial\hat{T}_t:=\{u\in\hat{T}_{\infty}\dvtx \mbox{ if }
u=u_1\cdots u_n \mbox{ for some }n,\\
&&{}\hspace*{71pt}\mbox{then }
u_1\cdots u_{n-1} \in\hat{T}_t, u\notin\hat{T}_t\}
\end{eqnarray*}
is the set of all external nodes or leaves in $\hat{T}_t$.

For this purpose we will study a closely related branching Markov
process (a~jump or step Markov process) $(T_t)_{t \geq0}$ which will
be described next.
We start the tree $T_0$ with one particle, the root \o, which is alive
at time $t=0$ and set accordingly $T_0:=\{\o\}$. This initial
ancestor dies at a random time $\tau_1$ where $\tau_1$ is exponentially
distributed with mean one and bears $b$ children. These children behave
independently from and similarly to their ancestor. After the first
birth the individuals $1,\ldots,b$ are alive with independent lifetimes
equal to\vspace*{-1pt} $(E_1^{(\o)},\ldots, E_b^{({\o})})$, respectively.
At time $\tau_2:=\min\{E_1^{(\o)},\ldots, E_b^{(
{\o})}\}+\tau_1$ the corresponding individual deceases and
gives birth to $b$ new individuals, namely its children. Because of the
memoryless property of the exponential distribution, all individuals
alive just after $\tau_2$ (namely the $b-1$ remaining children of the
root and the $b$ children of the root's child that deceased at time
$\tau_2$) behave similarly to and independently from each other, having
an exponentially distributed lifetime. We define $T_t$ as the tree
corresponding to the process described above when it is stopped at time $t>0$.
It is clear that we have $\hat{T}_{t}=T_{t+\tau_1}$.
In general, let $(\tau_j)_{0\leq j \leq\infty}$ be a sequence of
Markov times with
%
\begin{equation}\label{N_t}
\tau_0=0, \qquad \tau_j:=\min\{t\dvtx  N_t:=(b-1)j+1\},\qquad  N_t=|\partial
T_t|,
\end{equation}
where $N_t$ is the number of external nodes or leaves in the tree
$T_t$. Consequently, $\tau_j$ is the time of the $j$th death. By
convention, we consider all internal nodes at time~$t$ deceased and the
remaining external nodes alive. There are $k$ internal nodes if and
only if there are exactly $(b-1)k+1$ external ones.

Since $\tau_1 \stackrel{d}= \operatorname{Exp}(1)$, we obtain
\[
\tau_k - \tau_{k-1} \stackrel{d}= \min\bigl\{E_1, \ldots,
E_{(b-1)(k-1)+1}\bigr\}
\stackrel{d}=\operatorname{Exp}\bigl((b-1)(k-1)+1\bigr),
\]
where $(E_j)$ are identically distributed, independent random
variables, distributed as $E$ and representing the remaining lifetimes
of the $(b-1)(k-1)+1$ nodes alive
at time $\tau_{k-1}$ after the $(k-1)$th death.
For the distribution of $\tau_k$ we obtain then
%
\begin{equation}\label{tau}
\tau_k \stackrel{d}= \sum_{j=1}^k \frac{E_j}{(b-1)(j-1)+1}.
\end{equation}
\end{Remark}
\subsection*{General notation}
We denote for functions $f,g\dvtx\mathbb{N} \mapsto\mathbb{R}$ where
$g(n)\not=0\break\forall n \in\mathbb{N}$,
\[
f\sim_{\mathrm{a.s.}}g\quad  \Leftrightarrow\quad \lim_{n\rightarrow\infty}\frac
{f(n)}{g(n)}=1\quad \Leftrightarrow\quad  f(n)=g(n)\bigl(1+o(1)\bigr)
\]
with $o(1)\rightarrow0, n\rightarrow\infty$.
In the same manner we define
\[
f\leq_{\mathrm{a.s.}}g \quad \Leftrightarrow\quad \lim_{n\rightarrow\infty}\frac
{f(n)}{g(n)}\leq1\quad \Leftrightarrow\quad  f(n)\leq g(n)\bigl(1+o(1)\bigr)
\]
with $o(1)\rightarrow0, n\rightarrow\infty$.

We let $\mathbb{N}_0:=\mathbb{N}\cup\{0\}$.

We summarize the results in the following lemma which is the result of
\citeauthor{ChauvinKleinMarckertRouault05}
[(\citeyear{ChauvinKleinMarckertRouault05}), Lemma 2.1]
in the Yule tree case.
\begin{lemma}\label{lem2.3.1}
\textup{(1)} We have
$(\tau_k - \tau_{k-1})_{k\geq1}$ independent and
\[
\tau_k - \tau_{k-1} \stackrel{d}=\operatorname{Exp}\bigl((b-1)(k-1)+1\bigr).
\]

 (2) Further,
$(\tau_k)_{k\geq1}$ and $(T_{\tau_k})_{k\geq1}$ are independent.

 (3) We have
$(b-1)\tau_n\sim_{\mathit{a.s.}}\log n $ almost surely.
\end{lemma}
\begin{pf}
The first two statements can be proven with the same arguments given in
\citet{ChauvinKleinMarckertRouault05}. The last assertion is
Proposition 1 in Broutin and Devroye~(\citeyear{BroutinDevroye05}).
\end{pf}
%

Exactly in the same manner as for the tree $\hat{T}_\infty$, we assign
to each node in the tree~$T_t$ a label
$u \in\mathcal{U}$ in the natural way [cf. \refeq{eqlabels}].
Additionally, we weighted the edges in the tree $(\hat T_t)_{t\geq0}$
as the values induced by $\hat{T}_\infty$.
We define for $u= u_1\cdots u_n\in\mathcal{U},$
\[
D_u:=Z_{u_1}^{({\o})}+\sum
_{k=1}^{n-1}Z_{u_{k+1}}^{(u_1\cdots u_k)}
\]
as its weight or weighted position and accordingly for $l\in\mathbb
{Z}^d$ and $t\geq0$,
\[
\rho_t(l):=|\{u\in\partial T_t\dvtx D_u=l\}|.
\]
Note that we obtain for $t\geq0$, $\partial T_{t+\tau_1}$ is equal to
$\partial\hat{T}_t$.
%

%
\section{Continuous time martingales}\label{section6.2}
Our ambition is the study of the profile of a class of discrete time
trees $(\tilde{\mathcal{T}}_k)_k$ that can be constructed from the
class of continuous time trees considered in Section \ref{section6.1}
[cf. Section \ref{section6.6}]. We will follow the ideas used by \citet
{ChauvinKleinMarckertRouault05} in the case of the random binary
search tree and the Yule tree process
and construct trees $(T_t)_{t\geq0}$ described in the last section with
$T_{\tau_k}\stackrel{d}=\tilde{\mathcal{T}}_k$. Alternatively, we
deal with trees ${\mathcal{T}_k}:=T_{\tau_k}$
whose profiles are at least comparable to the profile of $\tilde
{\mathcal{T}}_k$. Recall that if we stop the tree process $T_t$ at time
$t=\tau_k$, we have $k$ internal nodes and $(b-1)k +1$ external nodes
or leaves in the tree.

These general trees are not as well studied as the binary search tree
and its continuous time analog, the Yule tree process.
\citet{ChauvinKleinMarckertRouault05}
used classical results from the Yule tree process and a corresponding
fragmentation process [cf.~Chauvin et al. (\citeyear{ChauvinKleinMarckertRouault05}), Section
2.2] to formulate their main profile convergence result
(see their Theorem 4.1).
In general, it seems difficult to use a connection with a fragmentation
processes.

Note that in the following we consider the natural filtrations
$(\mathcal{F}_t)_{t \geq0}$ and $(\mathcal{F}_{(k)})_{k\in\mathbb
{N}_0}$ where
\[
\mathcal{F}_t :=\sigma(T_s, s \leq t) \quad \mbox{and}\quad  \mathcal
{F}_{(k)}:=\sigma(T_{\tau_j}, 1 \leq j \leq k).
\]
For $v=(v_1,\ldots,v_d)$ and $w=(w_1,\ldots,w_d) \in\mathbb{C}^d$ we define
\[
v+w=(v_1+w_1,\ldots,v_d+w_d)  \quad \mbox{and}\quad    v\cdot w:=\sum
_{j=1}^dv_jw_j.
\]
%

%
Set
%
\begin{equation}\label{Z^t}
\tilde{Z}^{(t)}(dx):=\sum_{u\in\partial T_t}\bar\delta_{\{D_u\}}(dx),
\end{equation}
where $\bar\delta$ denotes the Dirac measure.
%

Let $m(\lambda)^t:= E\int e^{-\lambda\cdot x}\tilde{Z}^{(t)}(dx)$ for
$t \geq0$. Then, for $|m(\lambda)|<\infty$ we define the classical
martingale $W^{(t)}(\lambda)$ for $\lambda\in\mathbb{C}^d$, $t\geq0$, as
%
\begin{equation}\label{Wobent}
W^{(t)}(\lambda):=\frac{1}{m(\lambda)^{t}}\int e^{-\lambda\cdot
x}\tilde{Z}^{(t)}(dx).
\end{equation}
Using Biggins [(\citeyear{Biggins92}), page 148], it follows for $t\geq0$ that
\[
m(\lambda)^t=\exp[t(b Ee^{-\lambda\cdot Z}-1)].
\]
For the special case $d=1$ we define, setting $z:=e^{\lambda}\in\mathbb{C}$,
%
\begin{equation}\label{M(t,z)}
M(t, z):=W^{(t)}(-\lambda)=W^{(t)}(-\log z)= \sum_{u \in\partial
T_t}z^{D_u}e^{-t(b E z ^Z-1)}.
\end{equation}
This continuous time martingale was studied for $Z=1 $ and $ b=2$ (the
Yule tree case) in \citet{ChauvinKleinMarckertRouault05}.
%
\section{Discrete time martingales}\label{section6.3}
Let $\mathcal{T}_n:= T_{\tau_n}$ be the discrete time edge-weighted
tree with $n$ internal nodes. We define for $n \in\mathbb{N}$ and
$\lambda\in\mathbb{C}^d \setminus\{\lambda\dvtx \break C_n(\lambda)=0\}$,
%
\begin{equation}\label{W_n}
W_n(\lambda):= \frac{1}{C_n(\lambda)}\sum_{l \in\mathbb
{Z}^d}U_l(n)e^{-\lambda\cdot l}= \frac{1}{C_n(\lambda)}\sum_{u \in
\partial\mathcal{T}_n}e^{-\lambda\cdot{D_u}},
\end{equation}
where $C_n(\lambda)$ is a multiplicative factor which we will specify
below, and $U_l(n)$ is the number of external nodes in $\mathcal{T}_n$
at (weighted) level $l \in\mathbb{Z}^d$;
%
\begin{equation}\label{eqU_l(n)}
U_l(n):=\rho_{\tau_n}(l)=|\{u \in\partial T_{\tau_n} \dvtx D_u = l\}|,\qquad l\in\mathbb{Z}^d,\mbox{ } n\in\mathbb{N}.
\end{equation}
We choose $C_n(\lambda)$ in order to make $W_n(\lambda)$ a martingale
with respect to the filtration
$\mathcal{F}_{(n)}=\sigma\{\mathcal{T}_j, 1\leq j \leq n\}$.
In the case $d=1, b=2$ and $Z=1$, the random binary tree case, the martingale
%
\begin{equation}\label{Nummer}
M_n(z):=W_n(-\log(z)),
\end{equation}
where $z\in\mathbb{C}$ was found and $C_n(-\log z)$ calculated by
\citet{Jabbour-Hattab01}.

For the general case, let $D_n$ be the weighted depth of the $n$th
inserted (internal) node in $(\mathcal{T}_m)_{m\in\mathbb{N}_0}$
respectively $(T_t)_{t \in\mathbb{R}+}$.
Because of the exponential distribution of the lifetimes, each alive
individual (external node) is equally likely to be the next one to die
and to become an internal node.
Therefore we obtain for $l \in\mathbb{Z}^d$
\[
P(D_{n+1}=l|T_{\tau_n})=\frac{U_l(n)}{(b-1)n+1}.
\]
Assume now that node $u \in T_{\tau_n}$ is the next node to expire and
to bear $b$ new individuals in the tree
$T_{\tau_{n+1}}$. Denote by $Z_j^{(n+1)}, 1\leq j \leq b$, the
$d$-dimensional weight assigned to the edge from node $u$ to its $j$th child.
Then we have
%
\begin{equation}\label{eq4.1c}
U_l(n+1)-U_l(n)=-\mathbh{1}_{\{D_{n+1}=l\}} + \sum_{j=1}^b \mathbh{1}_{\{
D_{n+1}+Z_j^{(n+1)}=l\}}.
\end{equation}
With
$\tilde{W}_n(\lambda):=\sum_{l \in\mathbb{Z}^d}U_l(n)\exp(-\lambda
\cdot l)$ and $\tilde{W}_0(\lambda):=1$,
we obtain from \refeq{eq4.1c}
\begin{eqnarray*}
E\bigl(\tilde{W}_{n+1}(\lambda)|\mathcal{F}_{(n)}\bigr)
&= & \sum_{l }e^{-\lambda \cdot l}E\bigl(U_l (n+1)|\mathcal{F}_{(n)}\bigr) \\
&= & \tilde{W}_n(\lambda)- \sum_{l }e^{-\lambda\cdot l} \frac
{U_l(n)}{(b-1)n+1}{}\\
& &{} + \sum_{j=1}^b \sum_{l}e^{-\lambda \cdot l} P
\bigl(D_{n+1}=l-Z_j^{(n+1)}|\mathcal{F}_{(n)}\bigr) \\
&= &\tilde{W}_n(\lambda) \frac{(b-1)n + b Ee^{-\lambda \cdot Z}}{(b-1)n+1}.
\end{eqnarray*}
Iterating this we have that
\[
(W_n(\lambda))_{n \in\mathbb{N}}=\biggl(\frac{\tilde
{W}_n(\lambda)}{\prod_{j=0}^{n-1}((b-1)j+bEe^{-\lambda\cdot
Z})/((b-1)j+1)}\biggr)_{n\in\mathbb{N}}
\]
is a $(\mathcal{F}_{(n)})_{n \in\mathbb{N}}$ adapted
martingale. Consequently,
we set it according to our notation for $n\geq1$
%
\begin{equation}\label{eqC_n}
C_n(\lambda):= \prod_{j=0}^{n-1}\frac{(b-1)j+bE\exp(-\lambda\cdot
Z)}{(b-1)j+1},
\end{equation}
and $C_0(\lambda):=1$.
We define for further calculations and, for reference, the following set:
%
\begin{eqnarray}
&&N_C
:=\{\lambda\in\mathbb{C}^d\dvtx C_n(\lambda)=0 \mbox{ for some }
n\}\label{eqN_C}
\nonumber
\\[-8pt]
\\[-8pt]
\nonumber
&&\hspace*{21pt}=\biggl\{\lambda\dvtx -\frac{b}{b-1}E\exp(-\lambda\cdot Z) \in
\mathbb{N}_0\biggr\}.
\end{eqnarray}
%
%

\section{Relationships between discrete and continuous
time martingales}\label{section6.4}
In this section we will study the relationship of the two martingales
$(W_n(\lambda))_{n\in\mathbb{N}_0}$ and $(W^{(t)}(\lambda))_{t\geq0}$
defined in Sections \ref{section6.2} and \ref{section6.3}.\vspace*{1pt}
With $C_n(\lambda), n\in\mathbb{N}_0,$ defined in~\refeq{eqC_n}, we
set for $\lambda\in\mathbb{C}^d\setminus N_C$ [cf. \refeq{eqN_C}],
\[
\mathcal{H}_n(\lambda):=C_n(\lambda)e^{\tau_n(1-bEe^{-\lambda\cdot
Z})},\qquad  n\geq0.
\]
We claim that $(\mathcal{H}_n(\lambda))_{n \in\mathbb{N}_0}$ is a
martingale with respect to the filtration $(\mathcal{F}_{\tau_n})_{n
\in\mathbb{N}_0}$
with expectation 1.
\begin{lemma}\label{lemma2.3.2}
Let $\lambda\in\mathbb{C}^d \setminus N_C$.
Then:
\begin{longlist}[(1)]
\item[(1)]
$W^{(\tau_n)}(\lambda)= \mathcal{H}_n(\lambda) W_n(\lambda), \mbox{ }  n
\in\mathbb{N}_0.$
\item[(2)]
$(\mathcal{H}_n(\lambda))_{n \in\mathbb{N}_0}$ is a martingale with
respect to the filtration $(\mathcal{F}_{\tau_n})_{n \in\mathbb{N}_0}$
with expectation 1.
\item[(3)]
$(\mathcal{H}_n(\lambda))_{n \in\mathbb{N}_0}$ and $(W_n(\lambda
))_{n\in\mathbb{N}_0}$ are independent.
\end{longlist}
\end{lemma}
\begin{pf}
The first statement is a direct consequence of the definition. The
other two statements follow from an application of Lemma \ref{lem2.3.1}.
\end{pf}
%
%
\subsection{Asymptotic behavior and further relationships}
In this subsection we are interested in the asymptotic behavior of
$C_n(\theta)$ for $\theta\in\mathbb{R}^d$.
\begin{lemma}\label{lemma2.3.3} For $\theta\in\mathbb{R}^d$, $b>1$,
and $n\rightarrow\infty$ we have
\[
C_n(\theta) \sim_{a.s.} n^{{1}/{(b-1)}(bEe^{-\theta\cdot Z}-1)}\frac
{\Gamma({1}/{(b-1)})}{\Gamma({(bEe^{-\theta\cdot
Z})}/{(b-1)})}.
\]
\end{lemma}
\begin{pf}
Let $\alpha\not=0$ and $\beta\in\mathbb{R}$.
From Stirling's formula for the Gamma function [see, e.g., \citet
{FlajoletOdlyzko90}], we obtain
\begin{eqnarray*}
\prod_{j=0}^{n-1}\frac{\alpha j+\beta}{\alpha j+1}
& =& \pmatrix{n-(-\beta/\alpha) -1 \cr n}\left[\pmatrix{n-(-1/\alpha)-1
\cr
n}\right]^{-1}\\
& =& n^{({1}/{\alpha})(\beta-1)}\frac{\Gamma(1/\alpha)}
{\Gamma(\beta/\alpha)}\biggl(1+O\biggl(\frac{1}{n}\biggr)\biggr).
\end{eqnarray*}
With $\alpha=(b-1)$ and $\beta=bEe^{-\theta\cdot Z}$ we obtain the statement.
\end{pf}

Since $\mathcal{H}_n(\lambda)$ is a positive martingale for $
Ee^{-\lambda\cdot Z}>0$, we obtain immediately from a well-known,
classical martingale result, for every $\lambda\in\mathbb{R}^d$,
$\mathcal{H}_n(\lambda)$ converges almost surely to a limit $\mathcal
{H}(\lambda)$ as $n \rightarrow\infty$. More details are given in the
next lemma.\looseness=1
\begin{lemma}\label{lemma2.3.4}
Let $\lambda \in\mathbb{R}^d \setminus N_C$ with $N_C$ defined in
\refeq{eqN_C}. Then we have, almost surely,
\[
\mathcal{H}_n(\lambda) \rightarrow\biggl(\frac{Y}{b-1}\biggr)^{{(bE
e^{-\lambda\cdot Z}-1)}/{(b-1)}}
\frac{\Gamma({1}/{(b-1)})}{\Gamma({(bEe^{-\lambda
\cdot Z})}/{(b-1)})},\qquad   n\rightarrow\infty,
\]
where $Y\stackrel{d}=\Gamma(\frac{1}{b-1},\frac{1}{b-1})$.
\end{lemma}
\begin{pf}
First we show that $e^{-(b-1)t}N_t$, defined in \refeq{N_t}, converges
almost surely and that the limit law is gamma distributed.
From Athreya and Ney [(\citeyear{AthreyaNey75}), Remark 1, page 109] we detect that
\[
Es^{N_t}=s e^{-t}\bigl[1-\bigl(1-e^{-(b-1)t}\bigr)s^{b-1}\bigr]^{-1/(b-1)}.
\]
As a consequence of setting $\alpha:=b-1>0$ 
and $u=-is$ for $s\in\mathbb{R}$, we obtain by standard calculations
%
\begin{equation}\label{zwisch}
Ee^{isN_te^{-\alpha t}}\rightarrow\biggl(\frac{1}{\alpha is +1}
\biggr)^{{1}/{\alpha}},\qquad  t \rightarrow\infty.
\end{equation}
The right-hand side of \refeq{zwisch} is the characteristic function of
a $\Gamma(1/\alpha,1/\alpha)$ distributed random variable with density
\[
f_{\Gamma}(x):=\frac{(1/\alpha)^{1/\alpha}}{\Gamma(1/\alpha
)}x^{(1/\alpha)-1}e^{-x/\alpha}\mathbh{1}_{(0,\infty)}(x).
\]
With Doob's limit law it is immediately verified that the continuous
parameter nonnegative martingale
$(N_t e^{-(b-1)t})_{t\geq0}$\vspace*{1pt} converges with probability one
to a finite limit $Y$ which is gamma distributed
with parameters $(\frac{1}{\alpha},\frac{1}{\alpha})$:
%
\begin{equation}\label{eq4.2c}
N_te^{-(b-1)t} \rightarrow Y,\qquad   t\rightarrow\infty, \mbox{a.s.}
\end{equation}
Consequently, we have for $n\rightarrow\infty$,
%
\begin{equation}\label{eq6.3}
(b-1)n e^{-(b-1)\tau_n} \rightarrow Y \qquad\mbox{a.s.},
\end{equation}
since $\tau_n$ converges toward infinity almost surely.
Therefore, we obtain almost surely with Lemma \ref{lemma2.3.3} that
for $n\rightarrow\infty$ and $\lambda\in\mathbb{R}^d\setminus N_C$,
\begin{eqnarray*}
\hspace*{38pt}\mathcal{H}_n(\lambda)
& \sim_{\mathrm{a.s.}}& \bigl(n^{{1}/{(b-1)}} e^{-\tau_n}\bigr)^{(bEe^{-\lambda
\cdot Z}-1)}\frac{\Gamma({1}/{(b-1)})}
{\Gamma({(bEe^{-\lambda\cdot Z})}/{(b-1)})}\\
&\rightarrow_{\phantom{as}}& \biggl(\frac{Y}{b-1}\biggr)^{(bEe^{-\lambda
\cdot Z}-1)/(b-1)}\frac{\Gamma({1}/{(b-1)})}
{\Gamma({(bEe^{-\lambda\cdot Z})}/{(b-1)})}.\hspace*{38pt}\qed
\end{eqnarray*}
%
\noqed\end{pf}
Note that for $\tilde{Z}^{(t)}$ defined in \refeq{Z^t}, we have
\[
\tilde{Z}^{(t)}(dx)=\sum_{r=1}^{N_t} \mathbh{1}_{\{z_r^{(t)}\}}(dx),
\]
where $z_r^{(t)}$ is the weighted position of the $r$th individual alive
at time~$t$. If the $r$th individual is equal to node u, then $z_r^{(t)}$
is equal to $D_u\in\mathbb{Z}^d$.
Note that the positions of the individuals alive just before $\tau_n$---the
time of the $n$th death---are
$\{z_r^{(\tau_{n-1})}\dvtx 1\leq r \leq(b-1)(n-1)+1\}$ as the
particles do not move during their lifetimes.

Before we state the next theorem, we formulate the following
preliminary proposition.
The following set will play an essential role in the sequel
(cf. Theorem \ref{theorem2.3.5}).
For $1<\gamma\leq2$, $\lambda=\theta+i\eta$, define $\Omega_{\gamma
}^1, \Omega^2_{\gamma}\subset\mathbb{C}^d$ by
\begin{eqnarray*}
\Omega_{\gamma}^1&:=&\operatorname{int}\{\lambda\dvtx bEe^{-\gamma
\theta\cdot Z}<\infty\},\\
\Omega_{\gamma}^2&:=&\operatorname{int}\biggl\{\lambda\dvtx \frac{m(\gamma
\theta)}{|m(\lambda)|^{\gamma}}<1\biggr\}.
\end{eqnarray*}
We set
%
\begin{equation}\label{Nr}
\Lambda:=\bigcup_{1<\gamma\leq2} \Omega_\gamma^1\cap\Omega_\gamma^2.
\end{equation}

%
\begin{propos}\label{prop6.4.5}
For $\theta\in\Lambda$, $\delta>0$, $n\geq1$, $\alpha\in\mathbb
{R}$ we have:
\begin{longlist}[(1)]
\item[(1)]
$EW^{(\tau_n)}(\theta)=(Em(\theta)^{-\tau_n})^{-1}$;
\item[(2)]
$Ee^{-\theta z_1^{(\tau_n)}}=\frac{1}{Em(\theta)^{-\tau_n}((b-1)n+1)}$;\vspace*{2pt}
\item[(3)]
$E(\delta^{\alpha})^{\tau_n}\sim_{a.s.}\frac{\Gamma({(\log
(1/\delta)^\alpha+1)}/{(b-1)})}{\Gamma({1}/{(b-1)}
)}n^{{-\log(1/\delta)^\alpha}/{(b-1)}}.$
\end{longlist}
\end{propos}
%
%
\begin{pf}
For the proof of the first statement we let $n \geq1$, setting $a_n=(b-1)n+1$,
\begin{eqnarray*}
\xi_n&:=&\xi_n(\theta):=EW^{(\tau_n)}(\theta)=E\int e^{-\theta\cdot
x} \tilde{Z}^{(\tau_n)}(dx)=a_n \underbrace{Ee^{-\theta \cdot
z_1^{(\tau_n)}}}_{=:c_n}.
\end{eqnarray*}
In particular, we designate $\xi_1=\log m(\theta)+1$.
Using induction over $n$, we find
\begin{eqnarray*}
\xi_n
& = & E\Biggl[\sum_{r=2}^{a_{n-1}}e^{-\theta\cdot z_r^{(\tau_{n-1})}}
+ \sum_{j=1}^b e^{-\theta \cdot(z_r^{(\tau_{n-1})}+Z_j
)}\Biggr]{}\\
& = & \prod_{j=1}^n\frac{(b-1)(j-1)+1+\log m(\theta)}{(b-1)(j-1)+1}.
\end{eqnarray*}
%
%
Using \refeq{tau} we obtain
%
%
\begin{eqnarray}\label{eq6.6}
E(\delta^{\alpha})^{\tau_{n}}
& = & E \exp\Biggl( \sum_{j=1}^{n}\frac{E_j}{(b-1)(j-1)+1} \log\delta
^{\alpha}\Biggr)
\nonumber
\\[-8pt]
\\[-8pt]
\nonumber
& = & \prod_{j=1}^{n} \frac{(b-1)(j-1)+1}{(b-1)(j-1)+1 - \log\delta
^{\alpha}}.
\end{eqnarray}
%
%
From this we can finally conclude that $\xi_n=\xi_n(\theta)=
(Em(\theta)^{-\tau_n})^{-1}$.

For the second statement we obtain immediately from (1) that
\[
Ee^{-\theta\cdot z_1^{(\tau_{n})}}=c_{n}=\frac{1}{a_{n}}\xi_{n}=\frac
{1}{a_{n}}(Em(\theta)^{-\tau_{n}})^{-1}.
\]

Finally for the last assertion, using the same asymptotic result as in
Lemma \ref{lemma2.3.3}, we obtain
\begin{eqnarray*}
\prod_{j=1}^ n \frac{aj+1}{aj+1+x}
\sim_{\mathrm{a.s.}} (x+1)\frac{\Gamma({(x+1)}/{a})}{\Gamma
({1}/{a})}(n+1)^{{-x}/{a}},
\end{eqnarray*}
and further, using \refeq{eq6.6},
\begin{eqnarray*}
E(\delta^{\alpha})^{\tau_{n}}
&\sim_{\mathrm{a.s.}}& \frac{\Gamma({(\log(1/\delta)^{\alpha
}+1)}/{(b-1)})}{\Gamma({1}/{(b-1)})}n^{{-\log
(1/\delta)^{\alpha}}/{(b-1)}}.
\end{eqnarray*}
\upqed\end{pf}
In the following theorem we describe the convergence behavior of the
discrete time martingale of Section \ref{section6.3} [see \refeq{W_n}] and
of the continuous time martingale of Section \ref{section6.2} [cf.
\refeq{Wobent}].
The first part is based on Biggins (\citeyear{Biggins92}), Theorem 6. The proof of
the second part is
based on the application of the first part of Theorem \ref
{theorem2.3.5} and the relationship between the two martingales (cf.
Lemma \ref{lemma2.3.2}).

For a set $A \subset\mathbb{C}^d$ let $\operatorname{int}(A)$ denote the set
of interior points of $A$.
\begin{theorem}\label{theorem2.3.5}
With $\Lambda$ defined in \refeq{Nr} we have:
\begin{longlist}[(1)]
\item[(1)]
as $t\rightarrow\infty$, $\{W^{(t)}(\lambda)\}$ converges, a.s. and in
$L^1$, uniformly on every compact subset $C$ of $\Lambda$;
\item[(2)]
as $n\rightarrow\infty$, $\{W_n(\lambda)\}$ converges, a.s. and in
$L^1$, uniformly on every compact subset $C$ of $\Lambda$.
\end{longlist}
The limits are denoted as $W^{(\infty)}(\lambda)$, respectively,
$W_{\infty}(\lambda)$.
\end{theorem}
\begin{pf}
The first part was proven by Biggins [(\citeyear{Biggins92}), Theorem 6], so we only
need to prove the second part.

Let $C \subset\Lambda$ be a compact subset of $\Lambda$.
Therefore,
\[
\lim_{N}\sup_{n\geq N} E \sup_{\lambda\in C}|W_n(\lambda)-W_N(\lambda)|=0
\]
which implies the uniform $L^1$ convergence. Additionally, due to the
fact that $(\sup_{\lambda\in C}|W_n(\lambda)-W_N(\lambda)|)_{n\geq N}$
is a submartingale, this implies
also the a.s. uniform convergence.
From Lemma \ref{lemma2.3.2} we have
\[
W_n(\lambda)-W_N(\lambda)=E\bigl(W^{(\tau_n)}(\lambda) - W^{(\tau
_N)}(\lambda) |\mathcal{F}_{(n)}\bigr).
\]
Now taking the supremum and expectations we further deduce
\[
E \sup_{\lambda\in C}|W_n(\lambda)-W_N(\lambda)|\leq E\Bigl(\sup
_{\lambda\in C}\bigl|W^{(\tau_n)}(\lambda) - W^{(\tau_N)}(\lambda) \bigr| \Bigr).
\]
Taking the supremum over $n\geq N$ we get
\begin{eqnarray*}
\sup_{n\geq N}E \sup_{\lambda\in C}|W_n(\lambda)-W_N(\lambda)|
&\leq&E\sup_{n\geq N}\Bigl(\sup_{\lambda\in C}\bigl|W^{(\tau_n)}(\lambda)
- W^{(\tau_N)}(\lambda) \bigr| \Bigr)\\
&\leq& E\Delta_N,
\end{eqnarray*}
where we set
\[
\Delta_n:=\sup_{T \geq\tau_n}\Bigl(\sup_{\lambda\in C}
\bigl|W^{(T)}(\lambda)-W^{(\tau_n)}(\lambda)\bigr|\Bigr).
\]
Since $W^{(t)}(\lambda)$ converges a.s. uniformly, we have a.s. that
$\lim_{n\rightarrow\infty}\Delta_n=0$.
With the triangle inequality we obtain that $\Delta_n \leq2 \Delta_0$.
If we show that $\Delta_0$ is integrable, we deduce $\lim_{n}E\Delta
_n=0$ and
the statement of the theorem using dominated convergence.

Let $M_{\lambda}(s):=\{s>0\dvtx  W^{(s)}(\lambda)\not= W^{(s-)}(\lambda)\}
$ for $\lambda\in\Lambda$ and $s>0$.
As it is shown in Bertoin and Rouault
[(\citeyear{BertoinRouault03}), proof of Proposition 3] it
is sufficient for proving the
integrability of $\Delta_0$ to show that for all $x \in\Lambda$ a
disk,
\[
D_x({\bar\rho}):=\{\lambda\in\mathbb{C}^d\dvtx  \|\lambda-x\|<{\bar\rho
}\},
\]
exists with
\[
\sup_{ \lambda\in D_x({\bar\rho})}E\biggl(\sum_{s \in M_{\lambda}(s)
}\bigl|W^{(s)}(\lambda)-W^{(s-)}(\lambda)\bigr|^q\biggr)<\infty
\]
for some $q \in(1,2].$ Note that the set $M_\lambda(x)$ is a.s. countable.
We observe that $W^{(t)}(\lambda)\not=W^{(t-)}(\lambda)$ if and only if
$t=\tau_n$ for some $n\in\mathbb{N},$ and for that reason we can
conclude that
\begin{eqnarray*}
&&E\biggl(\sum_{s\in M_{\lambda}(s)} \bigl|W^{(s)}(\lambda
)-W^{(s-)}(\lambda)\bigr|^q\biggr) \\
&&\qquad{} =  E \sum_{n=1}^{\infty}\bigl|W^{(\tau_n)}(\lambda)-W^{(\tau
_n-)}(\lambda)\bigr|^q\\
&&\qquad{} =  \sum_{n=1}^{\infty}E \biggl|m(\lambda)^{-\tau_n}\biggl[\int
e^{-\lambda\cdot x}\,d\tilde{Z}^{(\tau_n)}(x)-
\int e^{-\lambda\cdot x}\,d\tilde{Z}^{(\tau_{n-1})}(x) \biggr]\biggr|^q.
\end{eqnarray*}
Now set
$a_{n-1}:=(b-1)(n-1)+1$ for $n\in\mathbb{N}$. Then we derive the
following formula for $\tilde{Z}^{(\tau_n)}$:
%
\begin{equation}\label{eq6.4b}
\tilde{Z}^{(\tau_n)}\stackrel{d}=\sum_{r=2}^{a_{n-1}}\mathbh{1}_{\{
z_r^{(\tau_{n-1})}\}}+\sum_{j=1}^b
\mathbh{1}_{\{z_1^{(\tau_{n-1})} + Z_j\}},
\end{equation}
where $(Z_1,\ldots, Z_b)$ are independent of $\{z_r^{(\tau
_{n-1})}, 1\leq r \leq a_{n-1}\}$ and distributed as $(Z_1^{(
{\o})}, \ldots, Z_b^{({\o})})$ [cf.
\refeq{vector}].

Then, with Lemma \ref{lem2.3.1}, Jensen's inequality, independence and
Proposition \ref{prop6.4.5},
we have for $q \in(1,2]$, $F_{\tau,n}:=\sigma(\tau_j, 1\leq j\leq n)$
and $n\in\mathbb{N}_0$, that
%
%
\begin{eqnarray}\label{eq4.1b}
&& E\bigl| W^{(\tau_n)}(\lambda)- W^{(\tau_n-)}(\lambda)\bigr|^{q}
 \nonumber\\
&&\qquad =  E\Biggl(|m(\lambda)|^{-q \tau_n}E\Biggl[\Biggl|e^{-\lambda \cdot
z_1^{(\tau_{n-1})}}\Biggl(\sum_{j=1}^b
e^{-\lambda \cdot Z_j}-1\Biggr)\Biggr|^{q} \mbox{ }\Big|\mbox{ }   F_{\tau,n}
\Biggr] \Biggr)
\nonumber
\\[-8pt]
\\[-8pt]
\nonumber
&&\qquad {\leq} E|m(\lambda)|^{-q \tau_n} E\bigl|e^{-\lambda z_1^{(\tau
_{n-1})}}\bigr|^q 2^{q-1}\Biggl(E \Biggl(\sum_{j=1}^b|
e^{-\lambda \cdot Z_j}|\Biggr)^{q}+1\Biggr)\\
&&\qquad=  \frac{2^{q-1}}{a_{n-1}} \frac{E|m(\lambda)|^{-q \tau
_n}}{Em(\theta q)^{-\tau_{n-1}}}
\Biggl( E \Biggl(\sum_{j=1}^b e^{-\theta \cdot Z_j}\Biggr)^{q} +1
\Biggr).\nonumber
\end{eqnarray}
%
%
Recall the definition of $\Lambda$ as
\[
\Lambda=\bigcup_{1<\gamma\leq2}(\Omega_{\gamma}^1\cap\Omega_{\gamma}^2).
\]
For $x\in\Omega_{\gamma}^1\cap\Omega_{\gamma}^2$
for some $\gamma\in(1,2]$, we can choose $\bar\rho$ sufficiently
small so that $D_{\bar\rho}(x) \subset \Omega_{\gamma}^1\cap\Omega
_{\gamma}^2$. Then there exist some $\delta<1$ so that for all
$\lambda\in D_{\bar\rho}(x)$ we have
$\frac{m(\theta\gamma)^{1/\gamma}}{|m(\lambda)|}\leq\delta$ for some
$\gamma\in(1,2]$.
Using this we obtain with Proposition \ref{prop6.4.5},
\begin{eqnarray*}
{ E|m(\lambda)|^{-\gamma\tau_n}(Em(\theta\gamma)^{-\tau_{n-1}})^{-1}}
\leq_{\mathrm{a.s.}}  C(\gamma, \lambda, b)n^{{\log\delta^{\gamma}}/{(b-1)}}.
\end{eqnarray*}
%
%
With \refeq{eq4.1b} by setting $q=\gamma\in(1,2]$ we obtain
\begin{eqnarray*}
 &&E\bigl| W^{(\tau_n)}(\lambda)- W^{(\tau_n-)}(\lambda)
\bigr|^{\gamma}  \\
&&\qquad \leq_{\mathrm{a.s.}}  2^{\gamma-1} \frac{1}{a_{n-1}} \Biggl(E\Biggl(\sum_{j=1}^b
e^{-\theta\cdot Z_j}\Biggr)^{\gamma}+1\Biggr)
C(\gamma, \lambda, b)n^{{\log\delta^{\gamma}}/{(b-1)}}\\
&&\qquad \sim_{\mathrm{a.s.}} \tilde C(\gamma, \lambda, b) n^{-(1+\varepsilon)}
\end{eqnarray*}
for a suitable constant
$\tilde C(\gamma,\lambda, b) $
independent of $n$, and
$\varepsilon:=\frac{\log(1/\delta) ^{\gamma}}{b-1}$.
Since $\sum_{n}n^{-1-\varepsilon}<\infty$ and $\sup_{\lambda\in D_{\bar
\rho}(x)}\tilde C(\gamma,\lambda, b)<\infty$
for ${\bar\rho}$ sufficiently small,
we obtain the statement.
\end{pf}

\begin{Remark}
The proof of the second part of Theorem \ref{theorem2.3.5} is roughly
the same as the corresponding proof of \citeauthor{ChauvinKleinMarckertRouault05}
[(\citeyear{ChauvinKleinMarckertRouault05}), Theorem
3.1]
up to the point where the integrability of $\Delta_0$ is proven.
Their arguments can also be used in our more general case except of the
one concerning
the integrability of $\Delta_0=\sup_{T \geq0}\sup_{\lambda\in
C}|W^{(t)}(\lambda)-1|$.
Since the Yule tree process
is a special kind of a fragmentation process, the integrability of
$\Delta_0$ in \citet{ChauvinKleinMarckertRouault05} can be obtained
from a result of \citet{BertoinRouault05}.
To be more precise for fragmentation processes, the integrability can
be verified by an application of the compensation formula for Poisson
point processes applied to the Poissonian construction of the
fragmentation (see \citeauthor{BertoinRouault03} [(\citeyear{BertoinRouault03}), Proposition 3] and
\citeauthor{Bertoin03} [(\citeyear{Bertoin03}), Theorem 2] for more details).
Is seems difficult to use a similar argument in the general case.
\end{Remark}
\begin{Remark}\label{remark_Lambda}
Recall that we have already defined
\[
\Lambda= \bigcup_{1<\gamma\leq2} (\Omega_{\gamma}^1\cap\Omega
_{\gamma}^2)
\]
with
\begin{eqnarray*}
\Omega_\gamma^1&= &\operatorname{int}\{\lambda\dvtx bEe^{-\gamma\theta\cdot
Z}<\infty\}\quad  \mbox{and}\\
\Omega_\gamma^2&= &\operatorname{int}\biggl\{\lambda\dvtx \frac{m(\gamma\theta
)}{|m(\lambda)|^{\gamma}}<1\biggr\}
\end{eqnarray*}
in \refeq{Nr}.
Using Biggins [(\citeyear{Biggins92}), Theorem 6], we have that $\{W^{(t)}(\lambda)\}
$ converges uniformly on any compact subset of $\Lambda$, almost surely
and in mean, as $t\rightarrow\infty$ [cf. Theorem \ref
{theorem2.3.5}].

Now define
\[
\tilde{\Lambda}:=\Lambda\cap\mathbb{R}^d
\]
as the restriction of $\Lambda$ on $\mathbb{R}^d$.
Then we can rewrite $\tilde\Lambda$ and show that
%
\begin{eqnarray}\label{tildeLambda}
\tilde{\Lambda}=\bigcup_{1<\gamma\leq2}{\Omega}_{\gamma}^1\cap\tilde
{\Omega}^3
\end{eqnarray}
with
\begin{eqnarray*}
\tilde{\Omega}^3&:=&\biggl\{\theta\in\mathbb{R}^d\dvtx\theta \in\Omega
_0:  -\log m(\theta)<\frac{-\theta\cdot m'(\theta)}{m(\theta)}\biggr\}
\quad  \mbox{and}\\
\Omega_0&:=&\operatorname{int}\{\lambda\in\mathbb{C}^d\dvtx m(\operatorname{Re}(\lambda
))<\infty\}
\end{eqnarray*}
(cf. Biggins [(\citeyear{Biggins92}), page 141]).
Since
\begin{eqnarray*}
m(\theta) & = & \exp(b E e^{-\theta\cdot Z}-1)\quad   \mbox
{and}\\
m'(\theta) & = & m(\theta)(-bEZe^{-\theta\cdot Z}),
\end{eqnarray*}
we have for any $\theta\in\Omega_0\cap\mathbb{R}^d$,
\[
\theta\in\tilde{\Omega}^3
 \quad \Leftrightarrow \quad 1 -b E e^{-\theta\cdot Z} < b \theta\cdot
EZe^{-\theta\cdot Z}.
\]
\end{Remark}
%
\subsubsection{Subtree sizes}\label{subsubsubtreesizes}
The rest of Section \ref{section6.4} is devoted to a further
characterization of $W^{(\infty)}(\cdot)$, resp. $W_{\infty}(\cdot)$,
defined in Theorem \ref{theorem2.3.5} as solutions of fixed point
equations. For deriving these fixed point equations we will split the
original tree in the $b$ subtrees which are growing from the children
of the root. For that purpose
we will investigate the sizes of the subtrees which are growing from a
node in the tree.
For every $u \in\mathcal{U}$ let
\[
\tau^{(u)}:=\inf\{t\dvtx  u \in T_t\}
\]
be the time of the first appearance (birth) of node $u$ in the tree.
For $t>0$ set
\[
T_t^{(u)}:=\{v \in\mathcal{U}\dvtx  uv \in T_{t+\tau^{(u)}}\},
\]
the tree process growing from $u$. Further set
\[
N_t^{(u)}=\bigl|\partial T_t^{(u)}\bigr|  \quad \mbox{and}\quad
n_t^{(u)}:=N_{t-\tau^{(u)}}^{(u)},
\]
the number of leaves at time $t\geq\tau^{(u)}$ in the subtree growing
from node $u$.
Then, using the same arguments as in \refeq{eq4.2c}, we obtain
%
\begin{eqnarray}\label{eq6.6b}
\lim_{t\rightarrow\infty}e^{-(b-1)t}N_t^{(u)} &=&Y _u
 \quad \mbox{and} \\
\lim_{t \rightarrow\infty} e^{-t(b-1)}n_t^{(u)}&=&Y_u e^{-\tau
^{(u)}(b-1)},\label{eq6.7b}
\end{eqnarray}
where $Y_u$ is distributed as $Y$; that is, it is $\Gamma(\frac
{1}{b-1},\frac{1}{b-1})$ distributed.\vspace*{1pt}
If $u,v$ are not in the same line of descent, we can conclude that by
the branching property~$Y_u$ and~$Y_v$ are independent. Since for $t
\geq\tau^{(u)}$ we have
\[
n_t^{(u)}=n_t^{(u1)}+n_t^{(u2)}+\cdots+n_t^{(ub)},\qquad  \tau^{(u1)}=
\tau^{(u2)}= \cdots=\tau^{(ub)},
\]
we obtain
%
%
using \refeq{eq6.7b} 
\[
e^{-\tau^{(u1)}(b-1)}\sum_{j=1}^b Y_{uj}=e^{-\tau^{(u)}(b-1)}Y_u,
\]
and for that reason we further get
%
\begin{equation}\label{subtree_sizes}
\frac{n_t^{(uj)}}{n_t^{(u)}}\rightarrow\frac{Y_{uj}}{\sum_{j=1}^b
Y_{uj}}=:U^{(uj)}, \qquad 1 \leq j \leq b.
\end{equation}
Finally, with $\tau_1=\tau^{(1)}$, we have
%
\begin{eqnarray}
Y&:=&Y_{\o}=e^{-\tau_1(b-1)}(Y_1+\cdots+Y_b),\label
{eq4.5c}\\
Y_j&=&U^{(j)}Y e^{\tau_1(b-1)},\qquad  1 \leq j \leq b
\quad \mbox{and}\label{eq4.6c}\\
1&=& \sum_{j=1}^b U^{(j)}\label{eq4.7c}.
\end{eqnarray}
The distribution of the subtree sizes and their limit distributions can
now be further calculated using a generalized P\'{o}lya--Eggenberger urn model.
%
%
\subsubsection{Limit martingale equation and splitting formulas}
In the following statement we derive representations of $W^{(\infty
)}(\cdot)$ and $W_{\infty}(\cdot)$ as solutions of fixed point
equations. Furthermore, the first part of Theorem \ref{theo2.3.6}
emphasizes the close relationship of $W^{(\infty)}(\cdot)$ and
$W_{\infty}(\cdot)$.
\begin{theorem}\label{theo2.3.6}
Let us assume $\lambda\in\tilde{\Lambda}$ [cf. \refeq{tildeLambda}].
Then the following formulas hold:
\begin{longlist}[(1)]
\item[(1)]
limit martingale connection,
%
\begin{eqnarray}\label{limitmartingaleconnection}
W^{(\infty)}(\lambda) &=&\biggl(\frac{Y}{b-1}\biggr)^{{(b Ee^{-\lambda
\cdot{Z}}-1)}/{(b-1)}}
\nonumber
\\[-8pt]
\\[-8pt]
\nonumber
&&{}\times\frac{\Gamma({1}/{(b-1)})}
{\Gamma({(b E e^{-\lambda\cdot{Z}} )}/{(b-1)})}W_{\infty
}(\lambda)\qquad  \mbox{a.s.},
\end{eqnarray}
with $Y \stackrel{d}=\Gamma(\frac{1}{b-1},\frac{1}{b-1} );$
\item[(2)]
splitting formula,
\begin{longlist}[(a)]
\item[(a)]
for the continuous time process
%
\begin{equation}\label{eqFixpunktgl1}
W^{(\infty)}(\lambda)=\sum_{j=1}^b e^{-\lambda\cdot{Z_j}}e^{-\tau
_1(bE e^{-\lambda\cdot{Z}}-1)}W_j^{(\infty)}(\lambda),
\end{equation}
where $W_1^{(\infty)}(\lambda), \ldots, W_b^{(\infty)}(\lambda) $ are
independent, distributed as $W^{(\infty)}(\lambda)$ and independent
of $\tau_1$ and $(Z_1, \ldots, Z_b)$ where
$(Z_1, \ldots, Z_b)$
are the weights assigned to the edge from node $\o$ to its
children $1, \ldots, b \in\mathcal{U}$,
\item[(b)]
for the discrete time process,
%
\begin{equation}\label{eqFixpunktgl2}
W_{\infty}(\lambda)=\sum_{j=1}^b e^{-\lambda\cdot{Z_j}}
\bigl(U^{(j)}\bigr)^{{(bE e^{-\lambda\cdot{Z}}-1)}/{(b-1)}}W_{\infty,
(j)}(\lambda),
\end{equation}
where $W_{\infty, (1)}( \lambda), \ldots,W_{\infty, (b)}( \lambda) $
are independent, distributed as $W_{\infty}(\lambda)$ and independent
of $(U^{(j)})$ where $(U^{(j)})$ are defined in (\ref{subtree_sizes}).
\end{longlist}
\end{longlist}
\end{theorem}
\begin{pf}(1)
This is a consequence of Lemmas \ref{lemma2.3.2} and \ref{lemma2.3.4}.
\begin{longlist}[(2)]
\item[(2)](a) For $t >\tau_1$ we have
\[
W^{(t)}(\lambda) 
 =  \sum_{j=1}^b e^{-\lambda\cdot Z_j}e^{-\tau_1(bEe^{-\lambda
\cdot Z}-1)}W_{j}^{(t-\tau_1)}(\lambda),
\]
where for $j=1,\dots, b$ we set
\[
W_{j}^{(t)}(\lambda):= \sum_{u \in\partial T_t^{(j)}}e^{-\lambda\cdot
D_u} e^{-t(bEe^{-\lambda\cdot Z}-1)},
\]
and we let $(Z_1,\ldots,Z_b)$ be the weights assigned to the edges from
$\o$ to its children $1, \ldots, b$.
Now let $t \rightarrow\infty$. Then the assertion follows.
\begin{longlist}
\item[(b)]
Follows now from (1), (2)(a) and \refeq{eq4.6c}.
\end{longlist}
\end{longlist}
\upqed\end{pf}
%
%
\section{Profiles and large deviation results}\label{section6.5}
To prove our main theorem, Theorem~\ref{theo2.3.9}, we will need some
preliminary lemmas that will be stated at the beginning of this
section. The first one, Lemma \ref{lemma2.3.9}, is based on Biggins
[(\citeyear{Biggins92}), Lemma 5] in the discrete time case with nonlattice
weights $Z$.
The proof for the continuous time nonlattice version can be managed
with some additional arguments [see Biggins (\citeyear{Biggins92}), page 150].
For the lattice case note that the critical points in the proof of
Lemma \ref{lemma2.3.9} are the values $\eta\in\mathbb{R}^d$
with
\[
\biggl|\frac{m(\lambda) }{m(\theta)}\biggr|=\exp\bigl(bEe^{-\theta\cdot
Z}\bigl(\cos(\eta\cdot Z)-1\bigr)\bigr)=1,\qquad \theta\in\Omega_0.
\]
This can only occur if $Z$ is a lattice distribution or if $\eta=0$.
Let
%
\begin{equation}\label{eqN}
N:=\bigl\{l\in\mathbb{Z}^d\dvtx P(\{Z=l\})>0\bigr\}
\end{equation}
denote the support of $Z$.
For $\eta=(\eta_1,\ldots,\eta_d) \in\mathbb{R}^d$ we denote by
\[
|\eta|:=\max\{|\eta_j|, j=1,\ldots, d\}
\]
the maximum-norm.
%
Denote for $v_1,\ldots, v_r \in\mathbb{R}^d$ by
\[
\mathcal{E}\{v_1,\ldots, v_r\}:=\Biggl\{v\dvtx \exists\lambda=(\lambda_1,\ldots
,\lambda_r) \in\mathbb{R}^r\dvtx v=\sum_{j=1}^r \lambda_jv_j\Biggr\}\subset
\mathbb{R}^d
\]
the subspace generated from the vectors $v_1,\ldots,v_r$ and,\vspace*{1pt}
analogously, let for some subset $S\subset\mathbb{R}^d,$ $\mathcal
{E}(S)$ be the subspace in $\mathbb{R}^d$ generated from the vectors in
the set~$S$.
The next proposition is obvious and will be stated for further
reference in the paper.
\begin{propos}\label{prop4.6.2}
Let $Z$ be a lattice distribution with support $N$; define
\[
f_{\theta}(\eta):=b Ee^{-\theta\cdot Z}\bigl(\cos(\eta\cdot Z)-1\bigr),\qquad  \theta
\in\Omega_0, \eta\in\mathbb{R}^d
\]
and let $a\in\mathbb{R}^1_+$.
\begin{longlist}[(1)]
\item[(1)]
If the dimension of $\mathcal{E}(N)$ is equal to $d$, then there are
only finitely many roots $\eta_0=0,\eta_1,\ldots,\eta_m,$ of $f_{\theta
}(\cdot)$ in $[-a,a]^d$
independently for all $\theta\in\Omega_0$ where $ m=m(a)\in\mathbb{N}_0$.
\item[(2)]
If the dimension of $\mathcal{E}(N)$ is $r<d$,
then there are $\eta_{0,1},\ldots,\eta_{0,d-r} \in\mathbb{R}^d$
(linearly) independent vectors,
$\eta_0=0,\eta_1,\ldots,\eta_{m} \in[-a,a]^d$ where $m:=\break m(a) \in
\mathbb{N}_0$ so that
\[
L:=\{\eta\in[-a,a]^d\dvtx f_{\theta}(\eta)=0\}\subset\bigcup_{j=0}^m
C_j,
\]
where
\[
C_j:=\eta_j + \mathcal{E}\{\eta_{0,1},\ldots, \eta_{0,d-r}\},\qquad 0\leq
j\leq m.
\]
\end{longlist}
\end{propos}
\begin{lemma}\label{lemma2.3.9}
For every lattice-distribution $Z$ we have almost surely
%
\begin{equation}\label{eq3}
\lim_{t\rightarrow\infty}\sup_{\theta\in\tilde K} \int_{|\eta|\leq
\pi}
\sqrt{t}^d\bigl|W^{(t)}(\theta+ i \eta) - W^{(\infty)}(\theta) \bigr|
\biggl|\frac{m(\theta+i\eta)}{m(\theta)}\biggr|^t\, d\eta=0
\end{equation}
for every compact subset $\tilde K \subset\tilde{\Lambda}$.
\end{lemma}

\begin{pf}
Let $C$ and $K$ represent here and in the rest of the proof some
arbitrary generic constants in $\mathbb{R}_+$ whose values may differ.
Let $f_{\theta}(\cdot)$ be defined as in Proposition \ref{prop4.6.2}.

Using the notation of Proposition \ref{prop4.6.2}, we will first
assume that the dimension of $\mathcal{E}(N)$ is equal to $d$.
As in \citet{Biggins92} we divide the integral into two parts for
$|\eta|<\varepsilon$ and $\varepsilon\leq|\eta|\leq{\pi}$.

We will consider the case $|\eta|<\varepsilon$ first.
With the standard Taylor series estimation we have for small $\varepsilon$
and $|\eta|< \varepsilon$ with $f_{\theta}(0)=f_{\theta}'(0)=0$,
%
\begin{eqnarray}\label{eqHilf}
\sup_{\theta\in\tilde K}\biggl|\frac{m(\lambda)}{m(\theta)}\biggr|\leq
\exp(-C|\eta|^2)
\end{eqnarray}
for some constant $C>0$.
Now let
\[
\tilde K_{\varepsilon}:=\{\lambda\dvtx \operatorname{Re}(\lambda) \in\tilde K, |\operatorname{Im}(\lambda)|\leq\varepsilon\}.
\]
Then, with Theorem \ref{theorem2.3.5}, we have
\begin{eqnarray*}
&&\sup_{\theta\in\tilde K}\int_{|\eta|<\varepsilon} \sqrt{t}^d
\bigl|W^{(t)}(\lambda)-W^{(\infty)}(\theta)\bigr|\biggl|\frac{m(\lambda)}{m(\theta
)}\biggr|^t\, d\eta \\
&&\qquad\leq\sup_{\lambda\in\tilde K_{\varepsilon}}\bigl|W^{(t)}(\lambda
)-W^{(\infty)}(\theta)\bigr|\int_{|\eta|<\varepsilon\sqrt{t}}\exp(-C|\eta|^2)\,
d\eta\\
&&\qquad\rightarrow K \sup_{\lambda\in\tilde K_{\varepsilon}}\bigl|W^{(\infty
)}(\lambda)-W^{(\infty)}(\theta)\bigr|,\qquad    t\rightarrow\infty.
\end{eqnarray*}
The last expression can be made arbitrarily small by choosing $\varepsilon
$ sufficiently small.

Next, we consider
\[
\sup_{\theta\in\tilde K} \int_{\varepsilon\leq|\eta|\leq\pi}
\sqrt{t}^d\bigl|W^{(t)}(\theta+ i \eta) - W^{(\infty)}(\theta) \bigr|
\biggl|\frac{m(\theta+i\eta)}{m(\theta)}\biggr|^t\, d\eta.\nonumber
\]
Let $\{\eta_1,\ldots,\eta_{m(\pi)}\}$ be the roots of $f_{\theta}(\eta
)$ or equivalently the values for which $|\frac{m(\theta+i\eta
)}{m(\theta)}|=1$, $\theta\in\Omega_0$ and $0<|\eta|\leq\pi$.\vspace*{1pt}
This is the case if and only if $\eta\cdot l\in2\pi\mathbb{Z}$ for
all $l\in N$ and $0<|\eta|\leq\pi$.
Consequently, we have $\cos(\eta_j \cdot Z)=1$, \mbox{$\sin(\eta_j \cdot
Z)=0$}, $e^{-i\eta_j \cdot Z}=1$ and $m(\theta)=m(\theta+i\eta_j)$ a.s.
for $j=1,\ldots, m=m(\pi)$.
With the triangle inequality we obtain
\begin{eqnarray*}
&&\sup_{\theta\in\tilde K} \int_{|\eta-\eta_j|<\varepsilon}
\sqrt{t}^d\bigl|W^{(t)}(\theta+ i \eta) - W^{(\infty)}(\theta) \bigr|
\biggl|\frac{m(\theta+i\eta)}{m(\theta)}\biggr|^t\, d\eta \\
&&\qquad\leq \sup_{\theta\in\tilde K} \int_{|\eta-\eta
_j|<\varepsilon}\sqrt{t}^d\bigl|W^{(t)}(\theta+ i \eta) - W^{(t)}(\theta) \bigr| \biggl|\frac
{m(\theta+i\eta)}{m(\theta)}\biggr|^t\, d\eta\\
&&\qquad\quad{}+ \sup_{\theta\in\tilde K} \int_{|\eta-\eta_j|<\varepsilon}
\sqrt{t}^d\bigl|W^{(t)}(\theta ) - W^{(\infty)}(\theta) \bigr| \biggl|\frac
{m(\theta+i\eta)}{m(\theta)}\biggr|^t\, d\eta\\
&&\qquad=: I_1^{(j)}+ I_2^{(j)} .
\end{eqnarray*}
For the second integral, $I_2^{(j)}$,
let $ \tilde K^{(j)}_{\varepsilon}:=\{\lambda=\theta+i\eta\in\mathbb
{C}^d\dvtx \theta\in\tilde K, |\eta-\eta_j|<\varepsilon\}$.
Again with a Taylor series estimation we have for $\varepsilon$ small and
$|\eta-\eta_j|<\varepsilon$,
\[
\sup_{\theta\in\tilde K}\biggl|\frac{m(\lambda)}{m(\theta)}\biggr|\leq
\exp(-C|\eta-\eta_j|^2).
\]
Then, for every $\varepsilon>0$ we have
\begin{eqnarray*}
&&\sup_{\theta\in\tilde K} \int_{|\eta-\eta_j|<\varepsilon}
\sqrt{t}^d\bigl|W^{(t)}(\theta ) - W^{(\infty)}(\theta) \bigr|\biggl |\frac
{m(\theta+i\eta)}{m(\theta)}\biggr|^t \,d\eta \\
&&\qquad\leq \sup_{\theta\in\tilde K^{(j)}_{\varepsilon}}\bigl|W^{(t)}(\theta
)-W^{(\infty)}(\theta)\bigr| \\
&&\qquad\quad{} \times\int_{|\eta-\eta_j|<\varepsilon}t^{d/2}\exp(-Ct|\eta
-\eta_j|^2)\, d\eta\\
&&\qquad= \sup_{\theta\in\tilde K}\bigl|W^{(t)}(\theta)-W^{(\infty)}(\theta
)\bigr|\underbrace{\int_{|\eta|<\varepsilon\sqrt{t}}\exp(-C|\eta|^2)\, d\eta
}_{\rightarrow(\pi/C)^{d/2}<\infty, t\rightarrow\infty}\rightarrow0,\qquad  t\rightarrow\infty,
\end{eqnarray*}
by Theorem \ref{theorem2.3.5}.

Next, we claim that $W^{(t)}(\theta+i\eta_j)=W^{(t)}(\theta)$ almost
surely for $j=1,\ldots, m$. For the proof note that $W^{(t)}(\theta
+i\eta_j)=W^{(t)}(\theta)$ because of $\eta_j \cdot Z_e\in2\pi\mathbb
{Z} $ a.s. for $e \in\pi(u), u\in\partial T_t $.
In particular we have also $W^{(\infty)}(\theta+i\eta_j)=W^{(\infty
)}(\theta)$ for $j=1,\ldots, m$.
Then, calculating the first integral, $I_1^{(j)}$, using \refeq{eqHilf},
\begin{eqnarray*}
&&\sup_{\theta\in\tilde K} \int_{|\eta-\eta_j|<\varepsilon}
\sqrt{t}^d\bigl|W^{(t)}(\theta+ i \eta) - W^{(t)}(\theta) \bigr| \biggl|\frac
{m(\theta+i\eta)}{m(\theta)}\biggr|^t\, d\eta \\
&&\qquad\leq\sup_{\lambda\in\tilde K^{(j)}_{\varepsilon}}\bigl|W^{(t)}(\lambda
)-W^{(t)}(\theta)\bigr| \int_{|\eta|<\varepsilon\sqrt{t}}\exp(-C|\eta|^2)\,d\eta\\
&&\qquad\rightarrow (\pi/C)^{d/2} \sup_{\lambda\in\tilde K^{(j)}_{\varepsilon
}}\bigl|W^{(\infty)}(\lambda)-W^{(\infty)}(\theta)\bigr|
\end{eqnarray*}
with Theorem \ref{theorem2.3.5} letting $t\rightarrow\infty$. This
can be made arbitrarily small by letting $\varepsilon\searrow0$ and,
therefore, $\eta\rightarrow\eta_j$.

Next we show that for $B:=\{\eta\dvtx \varepsilon\leq|\eta|\leq\pi, |\eta
-\eta_j|\geq\varepsilon, j=1,\ldots,m\}$, we obtain\vspace*{-1pt}
\[
\int_{B} \sqrt{t} ^d\bigl|W^{(t)}(\lambda)\bigr| \biggl|\frac{m(\lambda)}{m(\theta
)}\biggr|^t\, d\eta\rightarrow0, \qquad t\rightarrow\infty,\vspace*{-1pt}
\]
uniformly in a neighbourhood of any $\theta_0 \in\tilde{\Lambda}$
and hence uniformly on $\tilde K.$ The convergence of\vspace*{-1pt}
\[
\int_{B} \sqrt{t}^d \bigl|W^{(\infty)}(\theta)\bigr| \biggl|\frac{m(\lambda
)}{m(\theta)}\biggl|^t\, d\eta\vspace*{-1pt}
\]
toward zero will immediately follow from these considerations,
completing the proof.

 Now we have $B=\bigcup_{j\in I}\{\eta\dvtx a_j\leq|\eta|\leq b_j\}\subset
\{\varepsilon\leq|\eta|\leq\pi\}$ for some finite index set $I$ and
suitable $a_j< b_j$.
Additionally, $ \mu>1$ so that for all $j \in I$ we have $f_{\theta
}(\eta)\not=0$ if $\theta\in\Omega_0$ and $\frac{1}{\mu}a_j\leq|\eta
|\leq\mu b_j$.\vspace*{-1pt}
For $\theta_0 \in\tilde{\Lambda}$ let $B_{\bar\rho}:=\{\theta\in
\mathbb{R}^d\dvtx |\theta-\theta_0|\leq\bar\rho\}$ and let
$G^{(j)}_{c\bar\rho}:=\{\lambda\dvtx \theta\in B_{c\bar\rho}, c^{-1}
a_j \leq|\eta|\leq c b_j\}$, $j\in I$.\vspace*{-1pt}
It follows that $G^{(j)}_{c_1\bar\rho}\subset G^{(j)}_{c_2\bar\rho}$
for $0<c_1\leq c_2$, $j\in I$.
As $\theta_0\in\tilde{\Lambda}$, there is $\gamma\in(1,2]$ so that
$\theta_0 \in\Omega^2 _{\gamma}$. Then as a conclusion from the
definition of $\Omega^2_{\gamma}$ we have\vspace*{-1pt}
\[
\frac{ m(\gamma\theta_0)^{{1}/{\gamma}} }{ m(\theta_0) }<1.\vspace*{-1pt}
\]
Therefore, we can choose $\bar\rho$ sufficiently small, so that for
some $\delta<1$ we have $B_{ \mu\bar\rho}\subset\Omega_{\gamma}^2$,
and also\vspace*{-1pt}
\[
\frac{ \sup\{m(\gamma\theta)^{{1}/{\gamma}}\dvtx \theta\in B_{\mu
\bar\rho}\} }{\inf\{ m(\theta) \dvtx \theta\in B_{\mu\bar\rho}\}}\leq
\delta.\vspace*{-1pt}
\]
Further for $\theta\in\Omega_0$ we have $|\frac{m(\lambda
)}{m(\theta)}|<1$ for all $\frac{1}{\mu}a_j<|\eta| \leq\mu b_j$
and $j \in I$.\vspace*{1.5pt}
We obtain that for $r$ sufficiently small\vspace*{-1pt}
\[
\frac{\sup\{ |m(\lambda)|\dvtx \lambda\in G^{(j)}_{\mu r}\} }{\inf\{
m(\theta)\dvtx \theta\in B_{\mu r}\}}\leq\delta,\qquad j\in I.\vspace*{-1pt}
\]
Let $B^{(t)}(\lambda):=\int e^{-\lambda\cdot x} \tilde
{Z}^{(t)}(dx)-m(\lambda)^t$; then we obtain\vspace*{-1pt}
\[
W^{(t)}(\lambda) \biggl(\frac{m(\lambda)}{m(\theta)}\biggr)^t=m(\theta
)^{-t}B^{(t)}(\lambda) + \biggl(\frac{m(\lambda)}{m(\theta)}\biggr)^t.\vspace*{-1pt}
\]
Let $I_n:=\{t\dvtx n<t\leq n+1\}$. Then for $j\in I$ we find\vspace*{-1pt}
%
\begin{eqnarray}\label{eq4.3b}
&&\sup_{\theta\in B_{\bar\rho}, t\in I_n} \sqrt{t}^d \int
_{a_j\leq|\eta|\leq b_j} \bigl|W^{(t)}(\lambda) \bigr|\biggl|\frac
{m(\lambda)}{m(\theta)}\biggr|^t\, d\eta
\nonumber
\\[-8pt]
\\[-8pt]
\nonumber
&&\qquad\leq K\sqrt{n+1} ^d\biggl( \frac{\sup\{ |B^{(t)}(\lambda)|\dvtx \lambda
\in G^{(j)}_{\bar\rho}, t\in I_n \} }{\inf\{m(\theta)^{n+1}\dvtx \theta
\in B_{\bar\rho}\}}+\delta^n\biggr).
\end{eqnarray}
Now since $B^{(t)}(\lambda)$ is an analytic function, we can use
Cauchy's integral formula, the triangle inequality [see Biggins (\citeyear{Biggins92}), Lemma
3], and a compactness argument to
show that
%
\begin{eqnarray}\label{eq4.9b}
\sup_{\lambda\in G^{(j)}_{\bar\rho}}\bigl|B^{(t)}(\lambda)\bigr| \leq\frac
{1}{\pi} \sum_{\tilde{r}} \int_{C_{\tilde r}^{(j)}}
\bigl|B^{(t)}(z(s))\bigr|\,
ds,\qquad j\in I,
\end{eqnarray}
where $\{C_{\tilde{r}}^{(j)}\}$ parameterize the distinguished
boundaries of a finite number of disks covering $G^{(j)}_{\bar\rho}$
and lying within $G^{(j)}_{\mu{\bar\rho}}$ for $j\in I$.
As in \citet{Biggins92} we take expected values of \refeq{eq4.9b} and
use Jensen's inequality for some $\alpha> 1$ to obtain that for $j\in I$,
%
\begin{equation}\label{Ab}
  E \sup_{\lambda\in G^{(j)}_{\bar\rho}, t \in I_n}
\bigl|B^{(t)}(\lambda)\bigr| \leq K \sup_{\lambda\in G^{(j)}_{\mu{\bar\rho
}}}\Bigl( E\sup_{t\in I_n}\bigl|B^{(t)}(\lambda) \bigr|^{\alpha}\Bigr)^{{1}/{\alpha}}
\end{equation}
for some constant $K$.
Note that $m(\lambda)$ is never zero and that $|B^{(t)}(\lambda
)/m(\lambda)^t|$ is a regular submartingale. Hence with a
standard martingale inequality [Williams (\citeyear{Williams79}), Lemma 43.3], we obtain
\[
E \sup_{t \in I_n} \biggl|\frac{B^{(t)}(\lambda)}{m(\lambda)^t}
\biggr|^{\alpha} \leq\biggl(\frac{\alpha}{\alpha-1}\biggr)^{\alpha} E
\biggl|\frac{B^{(n+1)}(\lambda)}{m(\lambda)^{n+1}}\biggr|^{\alpha},
\]
and so the right-hand side of \refeq{Ab} is less than
%
\begin{eqnarray}\label{eq4.5b}
K \sup_{\lambda\in G^{(j)}_{\mu {\bar\rho}}}\bigl(
E\bigl|B^{(n+1)}(\lambda) \bigr|^{\alpha} \bigr)^{{1}/{\alpha}},\qquad j \in I.
\end{eqnarray}
Now the proof continues exactly as in Biggins [(\citeyear{Biggins92}),
Lemma 5],
bounding \refeq{eq4.5b} [cf. Biggins (\citeyear{Biggins92}), (4.6)]
with Biggins [(\citeyear{Biggins92}), Lemma 6] and finally shows that the expected
value of \refeq{eq4.3b} converges toward zero, and thus \refeq
{eq4.3b} converges toward zero almost surely, if $n\rightarrow\infty$.

For the second case when the dimension of $\mathcal{E}(N)$ is $r<d$ we
choose $a=\pi$ in Proposition \ref{prop4.6.2} and define
\[
D_j:=D_j(\varepsilon,\pi):=\{\eta\dvtx |\eta|\leq\pi,
d(\eta,C_j)<\varepsilon\},
\]
where
\[
d(\eta,C_j):=\min\{|\eta-w|, w\in C_j\}
\]
for $0\leq j\leq m=m(\pi)$ ($C_j$ was defined in the second part of
Proposition \ref{prop4.6.2}).
Note that we can choose $\delta>0$ so small that we have $m(\pi)=m(\pi
+\delta)$.
Set
\[
g_{\theta}(\eta)^{(t)}:=t^{d/2} \bigl|W^{(t)}(\theta+i\eta)\bigr|
\biggl|\frac{m(\theta+i\eta)}{m(\theta)}\biggr|^t.
\]

We have $|\frac{m(\theta+i\eta)}{m(\theta)}|<1$ for $\theta
\in B=B(\varepsilon,\pi):=[-\pi,\pi]^d \setminus\bigcup_{j=0}^m D_j$, and
\[
B(\varepsilon,\pi)\subset[-\pi-\delta,\pi+\delta]^d\Big\backslash\bigcup
_{j=0}^m D_j(\varepsilon/2,\pi+\delta)=:B(\varepsilon/2,\pi+\delta).
\]
We can now use similar arguments as in the case $r=d$ to obtain that
uniformly in a neighborhood of any
$\theta_0 \in\tilde{\Lambda}$, and hence uniformly on $\tilde K$ we have
\[
\int_{B} g_{\theta}(\eta)^{(t)} \, d\eta\rightarrow0,\qquad  t
\rightarrow\infty.
\]
\\
%
Let $\tilde{D}_j:=\{\lambda\in\mathbb{R}^{d-r}\dvtx \eta_j + \sum
_{k=1}^{d-r}\lambda_k \eta_{0,k} \in[-\pi,\pi]^d\}, 0\leq j\leq m$,
with\break
$(\eta_{0,k})_{1\leq k\leq d-r}$ as in Proposition \ref{prop4.6.2}.
Then we conclude that
\[
\sum_{j=0}^m\int_{D_j} g_{\theta}(\eta)^{(t)}\, d\eta\leq\sum_{j=0}^m
\int_{\lambda\in\tilde{D}_j}\int_{|\eta- ( \eta_j+\sum
_{k=1}^{d-r}\lambda_k \eta_{0,k})|<\varepsilon}
g_{\theta}(\eta)^{(t)}\,d\eta\, d\lambda.
\]
With the same calculations as in the case $r=d$, we get for every $\rho>0$,
\[
\int_{|\eta- ( \eta_j+\sum_{k=1}^{d-r}\lambda_k \eta_{0,k})|<\varepsilon}
g_{\theta}(\eta)^{(t)}\, d\eta< \rho
\]
if $\varepsilon\leq\varepsilon(\rho)$ and $t \geq T(\rho)$ independently of
the choice of $\eta_j,\lambda$ and $\theta$.
Therefore, for some $C\in(0,\infty)$ we have
\[
\sum_{j=0}^m\int_{D_j} g_{\theta}(\eta)^{(t)}\, d\eta\leq\sum_{j=0}^m
\int_{\lambda\in\tilde{D}_j} \rho\, d\lambda=C \rho,
\]
and this can be made arbitrarily small by choosing $\rho$ sufficiently
small. This completes the proof.
\end{pf}
\begin{lemma}\label{lemma2.3.10}
For any compact set $C \subset\tilde{\Lambda}$
we have almost surely
\[
\lim_{t\rightarrow\infty}\sup_{l \in\mathbb{Z}^d, \theta\in
C}e^{-\theta\cdot l} e^{t(1-b Ee^{-\theta\cdot Z})}t^{d/2}\bigl[ \rho
_t(l)-W^{(\infty)}(\theta)G_{(l,t)}(\theta)\bigr] = 0,
\]
with
\[
G_{(l,t)}(\theta):= \frac{1}{e^{-\theta\cdot l}}e^{t(bEe^{-\theta
\cdot Z}-1)}
\frac{1}{(2 \pi)^d}
\int_{|\eta|\leq\pi}e^{-btEe^{-\theta\cdot Z}(1-e^{i\eta\cdot
Z})}e^{-i\eta\cdot l}\,d\eta.
\]
\end{lemma}
\begin{pf}
Generally we have
$W^{(t)}(\theta)= e^{t(1-bEe^{-\theta\cdot Z})}\sum_l \rho
_t(l)e^{-\theta\cdot l}$ [cf. \refeq{Wobent}]
and
\begin{eqnarray*}
\biggl(\frac{1}{2\pi}\biggr)^d\int_{|\eta|\leq\pi}e^{i\eta\cdot
l}e^{-i\eta\cdot\tilde l}\,d\eta=
\cases{
0, &\quad  $\mbox{if }\tilde l\not= l,$\cr
1, &\quad  $\mbox{if }\tilde l=l.$}
\end{eqnarray*}
Using this and Lemma \ref{lemma2.3.9}, we obtain
\begin{eqnarray*}
&&\rho_t(l) e^{-\theta\cdot l} e^{t(1-bEe^{-\theta\cdot Z})}
\sqrt{t}^d\\
&&\qquad =  W^{(\infty)}(\theta) \sqrt{t}^d\frac{1}{(2\pi)^d}\int_{|\eta|\leq
\pi}e^{-btEe^{-\theta\cdot Z}(1-e^{i\eta\cdot Z})}e^{-i\eta\cdot
l}\,d\eta+ o(1),
\end{eqnarray*}
where the error term $o(1)$ is uniform in $l$ and in $\theta$ in any
compact subset of $\tilde\Lambda$.
\end{pf}
\begin{remark}
The next corollary deals with the special case when $d=1$ and when $Z$
takes only finitely many values in $\mathbb{N}_0$. In particular,
Corollary \ref{lemma2.3.11} gives detailed information about the term
$G_{l,t}(\theta)$ defined in Lemma \ref{lemma2.3.10}. Note that for
the proof of Corollary \ref{lemma2.3.11} we use the Cauchy formula to obtain
\[
\frac{1}{2\pi}\int_{-\pi}^\pi e^{btEe^{-\theta Z}e^{i\eta
Z}}e^{-i\eta l}\,d\eta=\frac{1}{l!}\frac{\partial^l}{\partial x^l}
\bigl(e^{btE(e^{-\theta}x)^Z}\bigr)_{|x=0}.
\]
Finally, standard calculations such as those used in Lemma \ref
{lemma2.3.10} lead then to the assertion of Corollary \ref{lemma2.3.11}.
\end{remark}
\begin{cor}\label{lemma2.3.11}
Let $d=1$ and $N\subset\{0,1,\ldots, L\}$ for some $L\in\mathbb{N}_0$. Then,
almost surely, for any compact $C \subset\tilde{\Lambda}$, we have
\[
\lim_{t\rightarrow\infty}\sup_{l \geq0,\theta\in C}e^{-\theta
l}e^{t(1-b Ee^{-\theta Z})}\sqrt{t}\biggl[ \rho_t(l)-W^{(\infty)}(\theta)
e^{-t} \frac{A_l}{l!} \biggr] = 0
\]
with
\[
A_l
 =  e^{btP(Z=0)}\sum_{D_l}\frac{l!}{a_1!a_2!\cdots a_L!}\prod_{j=1}^L
\bigl(btP(Z=j)\bigr)^{a_j},
\]
where $D_l:=\{(a_1, \ldots, a_L) \in\mathbb{N}_0^L\dvtx  \sum_{j=1}^L j
a_j =l\}.$
\end{cor}

 The terms in Corollary \ref{lemma2.3.11} can be further calculated.
Let $p_j:=P(Z=j)$ for $j \in\mathbb{N}_0$ and $\mathcal{P}^{(\xi
)}(l):=\frac{e^{-\xi}}{l!}\xi^l$ the Poisson measure with parameter
$\xi>0$.
Then
\begin{eqnarray*}
e^{-btEz^Z} \frac{z^l A_l}{l!}
& = & \sum_{\sum ja_j=l}e^{-bt(\sum_{j=1}^Lz^jp_j)}\prod_{j=1}^L \frac
{(btp_jz^j)^{a_j}}{a_j!}{}\\
& = & \sum_{\sum ja_j=l}\prod_{j=1}^L \mathcal{P}^{(tbz^jp_j)}(a_j).\\
\end{eqnarray*}
With the local limit theorem in Petrov [(\citeyear{Petrov75}), Theorem 7 of Chapter
VII, Section~2],
\[
\lim_{\lambda\rightarrow\infty}\sup_{l}\biggl|\sqrt{2\pi\lambda
}\mathcal{P}^{(\lambda)}(l)-\exp\biggl(-\frac{(l-\lambda)^2}{2\lambda
}\biggr)\biggr|=0,
\]
we have
\begin{eqnarray}\label{Ref}
\rho_t(l)
& = & z^{-l}e^{-t(1-bEz^Z)}t^{-1/2}\nonumber\\
&&{} \times\Biggl[ \sqrt{t}W(\infty,-\log(z))
\nonumber
\\[-8pt]
\\[-8pt]
\nonumber
&&\hspace*{16pt}{}\times\Biggl( \sum_{\sum ja_j=l}\prod_{j=1}^L \biggl\{\frac{1}{\sqrt{2\pi t b
z^jp_j}}\nonumber\\
&&{}\hspace*{88pt} \times\biggl(\exp\biggl\{-\frac
{(a_j-tbz^jp_j)^2}{2tbz^jp_j}\biggr\}+o(1)\biggr)\biggr\}\Biggr)+o(1)\Biggr].\nonumber
\end{eqnarray}
%
%
For this reason, using Corollary \ref{lemma2.3.11},\vspace*{1pt} we obtain an
expression which can be easily calculated for $A_l,l\geq1$ in the
binary search tree case [$Z=1$ a.s. and $A_l=(2t)^l$] [cf.
\citet{ChauvinKleinMarckertRouault05}], or also for the random
recursive tree case [$P(Z=0)=P(Z=1)=\frac{1}{2}$ and $A_l=e^t t^l$].

\citet{ChauvinKleinMarckertRouault05} calculated and
estimated the expression in \refeq{Ref} to prove their convergence result.
To prove our main theorem, Theorem \ref{theo2.3.9}, we chose
a~calculation inspired by a proof of \citet{Uchiyama82}.

Note that for\vspace*{1pt} $g\dvtx\mathbb{R}^d\mapsto\mathbb{R}$ we let $Dg(x)=g'(x)$
denote its gradient and $D^2g(x)$ its Hessian matrix at point $x \in
\mathbb{R}^d$, if such exists.
Further, for $M \in\mathbb{R}^d\times\mathbb{R}^d$ we denote by det
$M$ its determinant.
Let
%
\begin{equation}
A(-\theta):=bEe^{-\theta\cdot Z}-1= E\int e^{-\theta\cdot x} X(dx)-1
 \qquad\mbox{for }  \theta\in
\Omega_0
\end{equation}
with
%
\begin{equation}\label{eqX}
X:=\sum_{j=1}^b \mathbh{1}_{Z_j},
\end{equation}
where $(Z_1,\ldots,Z_b)\stackrel{d}=(Z^{(\o)}_1,\ldots,Z^{(\o)}_b$) are distributed as the
random weights attached to the edges that connect the root with its
children $1,\ldots,b$.
To prove our main Theorem \ref{theo2.3.9} we assume that $X$ is
nondegenerate in the sense that the support of the intensity measure
of $X$ is not contained in any $d-1$ dimensional hyperplane of $\mathbb{R}^d$.
Then $D^2A(-\theta)$ is a positive definite for all $\theta\in\Omega_0$,
by the inverse mapping theorem the set $\Omega_0^*:=\{DA(-\theta)\dvtx
\theta\in\Omega_0\}$
is open, and the mapping
\[
\theta\mapsto c=DA(-\theta)
\]
is a homeomorphism of $\Omega_0$ onto $\Omega_0^*$.
\begin{theorem}\label{theo2.3.9}
Let $K$ be a compact subset of $\Lambda^*:=\{DA(-\theta)\dvtx\theta\in
\tilde\Lambda\}\subset\mathbb{R}^d$ and assume that $X$ defined in
\refeq{eqX} is nondegenerate.
Then,
almost surely
\[
\lim_{n\rightarrow\infty}\sup_{c\in K}\biggl|\frac
{U_{l_n(c)}(n)}{A_c(n)}-W_{\infty}(\theta(c))\biggr|=0,
\]
where $\theta(c) \in\tilde{\Lambda}$
is chosen so that for
$c=:(c_1,\ldots,c_d)\in\mathbb{R}^d$
we have
\begin{eqnarray*}
DA(-\theta(c))&=&bEZe^{-\theta(c)\cdot Z}=c,\\
l_n(c)&:=&\biggl[\frac{c\log n}{b-1}\biggr]:=\biggl(\biggl[\frac{c_1\log
n}{b-1}\biggr],\ldots,\biggl[\frac{c_d\log n}{b-1}\biggr] \biggr)
\end{eqnarray*}
and
\begin{eqnarray*}
A_c(n)&:=&\frac{n^{{(bEe^{-\theta(c)\cdot Z}-1)}/{(b-1)}}}
{ e^{-\theta(c)\cdot{l_n(c)}}\sqrt{(2\pi{\log n }/{(b-1)}
)^d \det D^2A(-\theta(c))} }
\\
&&{}\times\frac{\Gamma({1}/{(b-1)})}{\Gamma({(bEe^{-\theta
(c)\cdot Z})}/{(b-1)})}.
\end{eqnarray*}
For the limit $W_{\infty}(\cdot)$ we have
%
\begin{equation}\label{fixedpoint}
W_{\infty}(\theta)=\sum_{j=1}^b e^{-\theta\cdot Z_j}
\bigl(U^{(j)}\bigr)^{{(bEe^{-\theta\cdot Z}-1)}/{(b-1)}}W_{\infty,
(j)}(\theta),
\end{equation}
where $W_{\infty, (1)}( \theta), \ldots,W_{\infty, (b)}( \theta) $ are
independent, distributed as $W_{\infty}(\theta)$ and independent of
$(U^{(j)})$ that are defined in (\ref{subtree_sizes}).
\end{theorem}
%
%
\begin{Remark}
From Theorem \ref{theorem2.3.5} it follows immediately from uniform
convergence on compact subsets that $(W_{\infty}(\lambda), \lambda\in
\Lambda)$ is a random analytic function.
Furthermore, we have for $ \theta\in\Lambda^*$, $W_{\infty}(\theta)$
is the unique solution of the fixed point equation \refeq{fixedpoint}
with expectation one. For this result, note that for $\theta\in\Lambda
^*$, $(W_n(\theta))_n$ is a nonnegative martingale with an (absolute)
first moment equal to one. From this we can conclude that $EW_{\infty
}(\theta)=1$ for all $\theta\in\Lambda^*$ (e.g., with Doob's limit law).
Finally, using the result of \citet{CaliebeRoesler05} we ascertain
that the solution of the fixed point equation \refeq{fixedpoint} with a
finite nonzero expectation is unique up to a~multiplicative constant.
\end{Remark}

We can also reformulate Theorem \ref{theo2.3.9} in terms of $l$
instead of in terms of $c$:
\begin{cor}\label{cor3.3.13}
Let $K$ be a compact subset of $\Lambda^*:=\{DA(-\theta)\dvtx\theta\in
\tilde\Lambda\}$.
Assume that $X$ defined in \refeq{eqX} is nondegenerate.
Then,
almost surely
\[
\lim_{n\rightarrow\infty}\sup_{l: {(b-1)}/{\log n}l\in K}
\biggl|\frac{U_{l}(n)}{\bar A_{{(b-1)l}/{\log n}}(n)}-W_{\infty}
(\theta_{l,n})\biggr|=0,
\]
where $\theta_{l,n} \in\tilde{\Lambda}$ is chosen so that
\[
bEZe^{-\theta_{l,n}\cdot Z}=\frac{(b-1)}{\log n}l\in\mathbb{R}^d
\]
and
\begin{eqnarray*}
\bar A_{{(b-1)l}/{\log n}}(n)
&:=&\frac{n^{{(bEe^{-\theta_{l,n}\cdot Z}-1)}/{(b-1)}}}{e^{-\theta
_{l,n}\cdot l} \sqrt{(2\pi{\log n}/{(b-1)})^{d}
\det D^2A(-\theta_{l,n})}}\\
&&{}\times\frac{\Gamma({1}/{(b-1)})}{\Gamma
({bEe^{-\theta_{l,n}\cdot Z}}/{(b-1)})}.
\end{eqnarray*}
Further, for $W_\infty(\cdot)$ we have
\[
W_{\infty}(\theta)=\sum_{j=1}^b e^{-\theta\cdot Z_j}
\bigl(U^{(j)}\bigr)^{{(bEe^{-\theta\cdot Z}-1)}/{(b-1)}}W_{\infty,
(j)}(\theta),
\]
where $W_{\infty, (1)}( \theta), \ldots,W_{\infty, (b)}( \theta) $ are
independent, distributed as $W_{\infty}(\theta)$ and independent of
$(U^{(j)})$ where $(U^{(j)})$ are defined in (\ref{subtree_sizes}).
\end{cor}

%
\begin{Remark}\label{rem6.5.9}
Note that the following two procedures are equivalent:
\begin{longlist}
\item[(1)]
take the supremum over $c \in K \subset\Lambda^*$ with $K$ a compact
subset and then choose $\theta(c)$, or
\item[(2)]
take the supremum over $\theta\in C \subset\tilde\Lambda$ with $C$ a
compact subset and then choose
$c(\theta)\dvt bEZe^{-\theta\cdot Z}=c(\theta)$.
\end{longlist}
\end{Remark}
\begin{pf*}{Proof of Theorem \protect\ref{theo2.3.9}}
For the proof we will use Lemma \ref{lemma2.3.10} and obtain
%
\begin{eqnarray}
&& e^{-\theta\cdot l}e^{t(1-bEe^{-\theta\cdot
Z})}t^{d/2}\rho_t(l) \label{eq4.19c}\nonumber\\
&&\qquad= W^{(\infty)}(\theta)
e^{t(1-bEe^{-\theta\cdot Z})}\biggl(\frac{\sqrt{t}}{2\pi}\biggr)^d\\
&&\qquad\quad{}\times\int
_{|\eta|\leq\pi} e^{t(bEe^{-\theta\cdot Z}e^{i\eta Z}-1)}e^{-i \eta
l}\,d\eta+o(1),\nonumber
\end{eqnarray}
where the error term is uniform for $\theta\in C \subset\tilde\Lambda$, a compact subset, and is uniform in $l$.
We claim that
\begin{longlist}[(a)]
\item[(a)]
\begin{eqnarray*}
&&e^{\tau_n(1-bE e^{-\theta(c)\cdot Z})}\biggl(\frac{\sqrt{\tau
_n}}{2\pi}\biggr)^d\int_{|\eta|\leq\pi} e^{\tau_n(bEe^{-\theta(c)\cdot
Z}e^{i\eta\cdot Z}-1)}e^{-i \eta\cdot l_n(c)}\,d\eta \\
&&\qquad =\frac{1}{\sqrt{(2\pi)^d \det D^2A(-\theta(c))}}+o(1),\qquad n\rightarrow
\infty,\nonumber
\end{eqnarray*}
with $o(1)$ uniform for $c$ in any compact subset $K \subset\Lambda
^*$, and
%
%
\item[(b)]\begin{eqnarray}\label{stetig}
&&\sup_{c \in K}\biggl|\biggl(\frac{Y}{b-1}\biggr)^{{(bEe^{-\theta
(c)\cdot Z}-1)}/{(b-1)}}
\nonumber
\\[-8pt]
\\[-8pt]
\nonumber
&&\quad{}\hspace*{6pt}\times\frac{\Gamma({1}/{(b-1)})}
{\Gamma
({(bEe^{-\theta(c)\cdot Z})}/{(b-1)})} W_{\infty}(\theta(c))
\biggr|<\infty
\end{eqnarray}
for every $K \subset\Lambda^*$ compact.
\end{longlist}
Since the functions in \refeq{stetig} are continuous in $\theta$ and
then also in $c$, (b) follows immediately.

The left-hand side of (a) is equal to
\begin{eqnarray*}
&&\tau_n^{d/2}e^{-\tau_n A(-\theta(c))}\frac{1}{(2\pi)^d}\int
_{|\eta|\leq\pi} e^{\tau_n A(-\theta(c) +i\eta)}e^{-i \eta\cdot
l_n(c)}\,d\eta\\
&&\qquad=\tau_n^{d/2}e^{-\tau_n
A(-\theta(c))}\bigl(I_1(\tau_n)+I_2(\tau_n)\bigr),
\end{eqnarray*}
where for $t\geq1$ we set
\begin{eqnarray*}
I_1(t)& = & \frac{1}{(2\pi)^d}\int_{\{|\eta|<\pi t^{-1/3}\}}e^{t
A(-\theta(c) +i\eta)}e^{-i \eta\cdot l_n(c)}\,d\eta \quad \mbox{and}\\
I_2(t)& = & \frac{1}{(2\pi)^d}\int_{\{\pi\geq|\eta|\geq\pi t^{-1/3}\}
}e^{t A(-\theta(c) +i\eta)}e^{-i \eta\cdot l_n(c)}\,d\eta.\\
\end{eqnarray*}
For $DA(-\theta(c))=bEZe^{-\theta(c)\cdot Z}=c$
we have with $\lambda=-\theta(c)+i\eta$
%
\begin{equation}\label{eq2.3}
A(\lambda) = A(-\theta(c))+ic\cdot\eta-\tfrac{1}{2}(D^2A(-\theta
(c))\eta)\cdot\eta+ i B |\eta|^3
\end{equation}
for $\eta\rightarrow0$ where $B=B(\lambda)$ is
uniformly bounded for $\eta\rightarrow0$ and $c \in K$.
With \refeq{eq2.3}, substituting $\eta=\mu/\sqrt{\tau_n}$, and using
\[
\int e^{-({1}/{2})(D^2A(-\theta(c))\mu)\cdot\mu}\, d\mu=( 2\pi)^d
( \operatorname{det} D^2 A(-\theta(c)) )^{-1/2},
\]
we obtain using Lemma \ref{lem2.3.1} and \refeq{eq2.3}, that
\begin{eqnarray*}
&&\tau_n^{d/2}e^{-\tau_n A(-\theta(c))}I_1(\tau_n )  \\
&&\qquad =  \frac{1}{(2\pi)^d}\int_{\{|\mu|<\pi\tau_n ^{1/6}\}}e^{i ({\mu}/{\sqrt{\tau_n }})\cdot(c\tau_n -[c\tau_n
+o(1)])}\\
&&\qquad\quad{}\times e^{-({1}/{2})(D^2A(-\theta(c))\mu)\cdot\mu} e^{iB|\mu|^3 \tau_n
^{-3/2}} \,d\mu\\
&&\qquad =  {} \frac{1}{\sqrt{(2\pi)^d \det D^2 A(-\theta(c) )}} +o(1),
\end{eqnarray*}
with $o(1)$ uniformly going to zero for $c \in K$.

Next we show $|t^{d/2}e^{-tA(-\theta(c))}I_2(t)|=o(1)$ uniformly for
all $c$ in a compact subset of $\Lambda^*$.
We have\vspace*{1pt}
\[
\bigl|t^{d/2}e^{-tA(-\theta(c))}I_2(t)\bigr|\leq t^{d/2}\int_{\{\pi\geq|\eta|\geq
\pi t^{-1/3}\}}e^{-q(\eta)t}\,d\eta\vspace*{1pt}
\]
with\vspace*{1pt}
\begin{eqnarray*}
q(\eta)&:= &A(-\theta(c))-\operatorname{Re}\bigl(A\bigl(-\theta(c) +i\eta\bigr)\bigr)\\[1pt]
&=&\tfrac{1}{2}
(D^2A(-\theta(c))\eta)\cdot\eta\bigl(1+o(1)\bigr),\qquad
\eta\rightarrow0,\vspace*{1pt}
\end{eqnarray*}
and $o(1)$ uniformly for all $c$ in a compact subset of $\Lambda^*$.
Additionally, if $\eta_1,\ldots,\eta_m$ are those values with $0< |\eta
|\leq\pi$ and $\eta\cdot Z \in2\pi\mathbb{Z}$ almost surely,
we have\vspace*{1pt}
\begin{eqnarray*}
q(\eta)&=& A(-\theta(c))-\operatorname{Re}\bigl(A\bigl(-\theta(c) +i\eta\bigr)\bigr)\\[1pt]
&=&bEe^{-\theta(c) \cdot Z}\bigl(1-\cos(\eta\cdot Z)\bigr)\\[1pt]
&=&bEe^{-\theta(c) \cdot Z}\tfrac{1}{2}|\eta-\eta_j|^2\bigl(1+o(1)\bigr), \qquad \eta
\rightarrow\eta_j, j=1,\ldots,m.\vspace*{1pt}
\end{eqnarray*}
Now note that $\operatorname{Re}(A(\lambda))<A(\theta(c))$ for $\pi\geq|\eta
|>\delta, \eta\not=\eta_j, j=1,\ldots,m$, for all $\delta>0$.
We can therefore choose $\varepsilon>0$ so small and independently of $c
\in K$ so that for all $\pi\geq|\eta|\geq\pi t^{-1/3} $ it follows that\vspace*{1pt}
\begin{eqnarray*}
q(\eta)t \geq
\min\biggl\{ \frac{\operatorname{det} D^2A(\theta(c))}{2}\pi^2 t^{1/3}, t \frac
{1}{2} bEe^{-\theta(c) \cdot Z}|\pi t^{-1/3}-\eta_j|^2,\tilde{C} t
\biggr\}
\geq\varepsilon t^{1/3}\vspace*{1pt}
\end{eqnarray*}
for some suitable constant $\tilde{C}>0$ and $t$ sufficiently large
(resp. $n$ if $t=\tau_n$).
It follows $|\sqrt{t}e^{-tA(-\theta(c))}I_2(t)|=o(1)$ and also $|\sqrt
{\tau_n}e^{-\tau_n A(-\theta(c))}I_2(\tau_n)|=o(1)$ with the error term
as claimed.
We obtain from
Theorem \ref{theo2.3.6} and \refeq{eq4.19c}, choosing $t=\tau_n$, $l=l_n(c)$,\vspace*{1pt}
\begin{eqnarray*}
&&e^{-\theta(c) \cdot l_n(c)}e^{\tau_n(1-bEe^{-\theta(c) \cdot
Z})}\tau_n^{d/2}\rho_{\tau_n}(l_n(c)) \\[1pt]
&&\qquad= \bigl(W_{\infty}( \theta(c) )+o(1)\bigr)\\[1pt]
&&\qquad\quad{}\times\biggl(
\frac{\Gamma({1}/{(b-1)})}{\Gamma({(bE e^{-\theta
(c)\cdot Z})}/{(b-1)})}\frac{({Y}/{(b-1)})^{
({bEe^{-\theta(c)\cdot Z}-1})/{(b-1)}}}{\sqrt{(2\pi)^d \det D^2A(-\theta
(c))}}\biggr).\vspace*{1pt}
\end{eqnarray*}
With \refeq{eq6.3} we have\vspace*{1pt}
\begin{eqnarray*}
&&e^{\tau_n(bEe^{-\theta(c)\cdot Z}-1)}\biggl(\frac{Y}{b-1}\biggr)^{
{(bEe^{-\theta(c)\cdot Z}-1)}/{(b-1)}}\\[1pt]
&&\qquad\sim_{\mathrm{a.s.}} n^{{(bEe^{-\theta(c)\cdot
Z}-1)}/{(b-1)}}\bigl(1+o(1)\bigr),\vspace*{1pt}
\end{eqnarray*}
and we obtain almost surely with $\tau_n\sim_{\mathrm{a.s.}}\frac{\log n}{b-1}, n\rightarrow\infty$,
\begin{eqnarray*}
&&\frac{e^{\tau_n(bEe^{-\theta(c)\cdot Z}-1)}(
{Y}/{(b-1)})^{{(bEe^{-\theta(c)\cdot Z}-1)}/{(b-1)}}}{e^{-\theta
(c)\cdot l_n(c)}\sqrt{(2\pi\tau_n)^d \det D^2A(-\theta(c))}}\frac
{\Gamma({1}/{(b-1)})}{\Gamma({(bEe^{-\theta
(c)\cdot Z})}/{(b-1)})} \\
&&\qquad \sim_{\mathrm{a.s.}}  \frac{n^{{(bEe^{-\theta(c)\cdot
Z}-1)}/{(b-1)}}}{e^{-\theta(c)\cdot l_n(c)}\sqrt{(2\pi{\log
n}/{(b-1)})^d \det D^2A(-\theta(c) )}}\\
&&\qquad\quad{}\times\frac{\Gamma(
{1}/{(b-1)})}{\Gamma({(bEe^{-\theta(c)\cdot Z})}/{(b-1)}
)}\\
&&\qquad =:  A_c(n).
\end{eqnarray*}
\upqed\end{pf*}
%
%
\begin{Remark}\label{rem3.3.15}
Consider the case $d=1$ and let $Z$ be bounded.
Define
\[
f(z)=1-bEz^Z+\log(z)bEZz^Z
\]
for $z>0$.
Then it follows immediately that $\tilde{\Lambda}=\{\theta\in\mathbb
{R}\dvtx  f(e^{-\theta})<0\}$ [cf. Remark~\ref{remark_Lambda}].
Since $f'(z)=\frac{1}{z}(\log(z)bEZ^2z^Z)=0\Leftrightarrow
z=1$ we have
a single local minimum ($f''(z)>0$) at the point $z=1$ with $f(1)=1-b$.
\end{Remark}

Further, for $Z\geq0$ a.s. we have
\[
\lim_{z \searrow0}f(z)=1-bp_0,\qquad  \lim_{z \rightarrow\infty
}f(z)=\infty.
\]
Consequently, we have one root of $f$ if $p_0>\frac{1}{b}$ and
otherwise two roots of $f$. In the first case, let $z_0=0$ and let
$z_1$ be the root of
$f$. In the second case, let $z_0<z_1$ be the two roots of $f$. In both
cases we have $\tilde\Lambda=(-\log(z_1),-\log(z_0))$ [where we set
$-\log(0):=\infty$].
%
\section{Examples}\label{section6.6}
Note that most of the following examples are taken from \citet
{BroutinDevroye05}.

In order to simplify notation and to work out various connections to
known results, in the case $d=1$ and for $z \in\mathbb{R}^+$, we use
$M_{\infty}( z):=W_{\infty}(-\log(z))$ instead of $W_{\infty}(\lambda)$
[cf. \refeq{Nummer}]. Set
\begin{eqnarray*}
V&:=&\{e^{-\lambda}\dvtx \lambda\in\Lambda\}, \\
\tilde{V}&:=& V \cap\mathbb{R}\quad  \mbox{and}\\
V^* &:=&\{bE(Zz^Z)\dvtx z\in\tilde{V}\}={\Lambda}^*.
\end{eqnarray*}
In complete analogy to Corollary \ref{cor3.3.13} we have the following:
\begin{cor}\label{cor4.6.8}
Let $d=1$ and let $K$ be a compact subset of
\[
V^*:=\{bE(Zz^Z)\dvtx z \in\tilde V\}.
\]
Then almost surely
\[
\lim_{n\rightarrow\infty}\sup_{l: {(b-1)l}/{\log n}\in K}
\biggl|\frac{U_{l}(n)}{\hat A_{{(b-1)l}/{\log n}}(n)}-M_{\infty}
(z_{l,n})\biggr|=0,
\]
where $z_{l,n} \in\tilde{V}$ is chosen so that
\begin{eqnarray*}
bEZz_{l,n}^Z&=&\frac{l(b-1)}{\log n},
\\
\hat A_{{(b-1)l}/{\log n}}(n)&:=&\frac{n^{
{(bEz_{l,n}^Z-1)}/{(b-1)}}}{z_{l,n}^{l}\sqrt{2\pi\log n {b}/{(b-1)}
E(Z^2z_{l,n}^Z)}}\frac{\Gamma({1}/{(b-1)})}{\Gamma
({bEz_{l,n}^Z}/{(b-1)})}
\end{eqnarray*}
and
\[
M_{\infty}(z)=\sum_{j=1}^b z^{Z_j}\bigl(U^{(j)}\bigr)^{
{(bEz^{Z}-1)}/{(b-1)}}M_{\infty, (j)}(z),
\]
where $M_{\infty, (1)}( z), \ldots,M_{\infty, (b)}( z) $ are
independent, distributed as $M_{\infty}(z)$ and independent of
$(U^{(j)})$ where $(U^{(j)})$ is defined in (\ref{subtree_sizes}).
\end{cor}
\begin{Example}[(Random binary search tree)]
A random binary search tree can be built incrementally. Let $U_1,\ldots
, U_n$ be independent random variables uniformly~distributed over the
unit interval.
We start the tree by storing $U_1$ in the root node. If $U_2$ is
greater than $U_1$, we add a right child to the root and store $U_2$ in
that node. If $U_2$ is less than $U_1$, we add a left child to the root
and store $U_2$ in that node. Then we repeat that procedure
incrementally for $U_3,\ldots, U_n$.
The nodes where we stored some $U_j$ for some $j$ are called internal nodes.
We refer to \citet{Devroye91}, \citet{Devroye98a} and the references
given there for the construction of binary search trees. A summary of
known results about binary search trees is given in \citet{Mahmoud92}
and \citet{Knuth98}.

Let $\mathcal{T}_n$ be a random binary search tree with $n$ (internal)
nodes. We will only consider complete binary search trees. That means
that we add $n+1$ external nodes to each binary search tree with $n$
internal nodes in the following manner.
If $u$ is an internal node and has no offspring, we add two external
nodes as its potential children to it. If it has already one child,
then we add one external node to $u$ as a second potential child. If
$u$ has already two children, we add nothing.
Note that every external node corresponds to one of the free places
available for the sorting of a new internal node and that each free
place is likely to be chosen next with equal probability.

It is well known that for the random binary search tree~$\mathcal{T}_n$
we have
\[
\lim_{n\rightarrow\infty}\frac{|\min\{D_u\dvtx u \in\partial\mathcal
{T}_n\}|}{\log n}=\alpha_-,
\qquad
\lim_{n\rightarrow\infty}\frac{|\max\{D_u\dvtx u \in\partial\mathcal
{T}_n\}|}{\log n}=\alpha_+,
\]
where $\alpha_-,\alpha_+$ are the only nonnegative solutions of the
equation $x\log\frac{x}{2}-x+\break 2=1$ [see, e.g., Devroye (\citeyear{Devroye86,Devroye87,Devroye98a}) and
references given there].
Chauvin et al. (\citeyear{ChauvinKleinMarckertRouault05}) proved the
following result, which is covered by Theorem \ref{theo2.3.9}:
\begin{theorem}[{[\citet{ChauvinKleinMarckertRouault05}]}]
Almost surely, for any compact subset~$K$ of $(\alpha_-,\alpha_+)$,
\[
\lim_{n\rightarrow\infty} \sup_{l:{l}/{\log n} \in K}\biggl(\frac
{U_l(n)}{EU_l(n)}-M_{\infty}\biggl(\frac{l}{2\log n}\biggr)\biggr)=0.
\]
\end{theorem}

The profile of the binary search tree was also studied with other methods.
Let $\alpha$ denote the limit ratio of the level and the logarithm of
the tree size.
Then,
Fuchs, Hwang and Neininger (\citeyear{FuchsHwangNeininger05}) proved
convergence in distribution for $\alpha\in V^*=(\alpha_-,\alpha_+)$
and for $\alpha\in[1,2]$ convergence of all moments to prove their results.
They used the contraction method and the method of moments.

Drmota, Svante and Neininger (\citeyear{DrmotaJansonNeininger06}) treated a
class of generalized $m$-ary search trees including binary search
trees. For those trees they proved that in a~certain range the
normalized profile converges in distribution (Theorem 1.1).
They used arguments based on the contraction method in order to prove
convergence in distribution of several random analytic functions in a
complex domain.
\end{Example}
\begin{Example}[(Random recursive tree)]\label{ex3.3.2}
A random recursive tree is built inductively. The tree $\tilde{\mathcal
{T}}_1$ consist of a single node \o, the root. Let $\tilde{\mathcal
{T}}_n$ already exist and consist of the nodes $\{v_1, \ldots, v_n\}
$. To grow the tree choose a node $v_j$ out of the set $\{v_1, \ldots,
v_n\}$ uniformly and at random, and attach the new node $v_{n+1}$ as a
child to node $v_j$ [cf. \citet{SmytheMahmoud95} and references given
there].

Fuchs, Hwang and Neininger (\citeyear{FuchsHwangNeininger05}) showed that
the profile of the random recursive tree
normalized by its mean converges in distribution if the limit ratio
$\alpha$ of the level and the logarithm of the
tree size lies in $[0,e)$. They also showed convergence of all moments
to hold for $\alpha\in[0,1]$. Furthermore, they proved that inside
the interval $(1,e)$ only convergence of a finite number of moments is possible.
Drmota and Hwang (\citeyear{DrmotaHwang05b}) showed that the variance of
the profile $U_l(n)$ of the random recursive tree asymptotically
undergoes four phase transitions and exhibits a bimodal behavior in
contrast to the unimodality of the expected value of the profiles (cf.
comments made on this topic in the \hyperref[intr]{Introduction}). For $l$ around the
most numerous level (where the width is attained), the value of the
martingale shall be a.s. constant; more precisely one has $M_{\infty
}(l/\log n)=1$ almost surely [cf. Drmota and Hwang (\citeyear{DrmotaHwang05b,DrmotaHwang05a})].
In the sequel, \citet{DrmotaHwang05a} sketched
that $U_l(n)\sim_{\mathrm{a.s.}} M_{\infty}(\alpha)EU_l(n)$ almost surely, where
$\alpha=\lim_{n}\frac{l}{\log n}\in[0,1)$, using a martingale argument
of \citet{ChauvinDrmotaJabbour01} and Cauchy's integral formula.
%

We will show below as an application of Theorem \ref{theo2.3.9} that
the profile of the random recursive tree
normalized by its mean converges almost surely if the limit ratio
$\alpha$ of the level and the logarithm of the
tree size lies in $(0,e)$. Additionally the profile converges uniformly
for $\alpha$ in any compact subset of $(0,e)$.

First note that it is possible to interpret a random recursive tree
with $n$ internal nodes as a weighted binary tree $T_{\tau_{n-1}}$ with
$n-1$\vspace*{-2pt} internal nodes by weighting the edges with independent copies of
$Z\stackrel{d}= \operatorname{Bernoulli} (\frac{1}{2})$ and finally
by interpreting the~$n$ external nodes of the latter as the internal
nodes of the former. We have to choose $Z_2=1-Z_1$.
This follows immediately since every external node in the weighted
binary tree is equally likely to be the next one to die and to get two
external children where in the recursive tree each internal node is
equally likely to be the next one to produce an offspring.
For more details on the construction we refer to Broutin and Devroye
(\citeyear{BroutinDevroye05}), Section 4.2.

With this construction it is clear that not only is the height treated
in \citet{BroutinDevroye05}, but also the
distribution of the profile is kept by this construction. By embedding
the random recursive tree $(\tilde{\mathcal{T}}_n)_{n\geq1}$ in the
weighted tree process and by identifying
$(T_{\tau_{n-1}})_{n\geq1}$ with $(\tilde{\mathcal{T}}_n)_{n\geq1}$,
we deduce the following:
\begin{theorem}
Let $K\subset(0,e)$ be a compact subset. Then almost surely
\[
\lim_{n\rightarrow\infty}\sup_{l:{l}/{\log n}\in K}\biggl|\frac
{U_{l}(n)}{EU_{l}(n)}-M_{\infty}\biggl(\frac{l}{\log n}\biggr)\biggr|=0,
\]
and for $z \in(0,e)$,
\[
M_{\infty}(z)=z U^{z}M_{\infty,(1)}(z)+(1-U)^{z}M_{\infty,(2)}(z),
\]
where $M_{\infty,(i)}(z)\stackrel{d}=M_{\infty}(z), U$ is
$\operatorname{uniform}[0,1]$ random variable and
$M_{\infty}(z)$, $M_{\infty,(1)}(z)$, $M_{\infty,(2)}(z)$ and $U$ are
independent.
\end{theorem}
\begin{pf}
Obviously, using Corollary \ref{cor4.6.8}, we have:
\begin{longlist}[(1)]
\item[(1)]
$\tilde{V}=\{z>0\dvtx 1-(1+z)<-\log(z)z\}=(0,e);$
\item[(2)]
$V^*=\{2 E(Zz^Z)=z \dvtx z \in\tilde V\}=(0,e);$
\end{longlist}
and we have $z_{l,n}= \frac{l}{\log n}$.
Then with
\[
{ \hat A}_{l/\log n}(n)
 =  \frac{n^{l/\log n}}{\Gamma(1+l/\log n)(l/\log n)^{l}\sqrt{2\pi l }}
\]
and
\[
B_l(n):=\frac{(\log n)^l}{\Gamma(1+{l}/{\log n})\sqrt{2\pi l}}
\biggl(\frac{e}{l}\biggr)^l,
\]
we easily obtain
\[
\frac{\hat{ A}_{l_n/\log n}(n)}{B_{l_n}(n)}=1+o(1).
\]
Finally, we note that Hwang (\citeyear{Hwang95}) showed that [see
also\break
\citeauthor{FuchsHwangNeininger05} (\citeyear{FuchsHwangNeininger05}), equation (3), page 2]
\[
EU_l(n)=\frac{(\log n) ^l}{l!\Gamma(1+{l}/{\log n})}\bigl(1+o(1)\bigr)\sim
_{\mathrm{a.s.}}B_l(n)\bigl(1+o(1)\bigr)
\]
which yields the theorem by using the Stirling formula.
%

From Theorem \ref{theo2.3.6}, part~(2)(b), using $Z_1=1-Z_2$,
$U^{(1)}+U^{(2)}=1$ and by setting $U:=U^{(1)}$, we obtain
\[
M_{\infty}(z)=z U^{z}M_{\infty,(1)}(z)+(1-U)^{z}M_{\infty,(2)}(z)
\]
and that the claimed independence relations also hold.
For the distribution of $U$ note that if $E_1,E_2$ are independent,
exponentially distributed random variables, then $\frac{E_1}{E_1+E_2}$\vspace*{1pt}
is $\operatorname{uniform}[0,1]$ distributed. Now $Y$, defined after Lemma \ref
{lemma2.3.3}, is Gamma distributed with parameters (1, 1) which is the
same as being exponentially distributed with expectation one.
It follows that $U^{(i)}, i=1,2$, is uniformly~$[0,1]$ distributed.
\end{pf}

 Note that $V^*=(0,e)$ is the natural range for convergence, since,
\citet{Devroye87} and \citet{Pittel94} showed for the height $H_n$ of
$\tilde{\mathcal{T}}_n$ that
\[
\lim_{n\rightarrow\infty}\frac{H_n}{\log n}\stackrel{P}{\rightarrow}e.
\]
So $e$ should be the upper bound for any range of convergence of the profile.
\end{Example}
\begin{Example}[(Random lopsided trees)]\label{ex3.3.20}
Prefix-free codes are particularly interesting because they can be
decoded directly by following a path in a tree and output a character
corresponding to the codeword when reaching a leaf. Each node~$u$
represents a prefix $p$ and its children represent the words that can
be built by appending a symbol to $p$. When reaching a leaf, one
obtains a character corresponding directly to the codeword.

Some codes have encoding length depending on the symbols. These codes
are called Varn codes [cf. \citet{Varn71}] and naturally lead to
lopsided trees. Lopsided trees are trees with edges having nonequal length.
We refer to \citet{BroutinDevroye05} for further details, especially
on the height of such trees and for further references.
There are no results about the asymptotic behavior of the profile of
random lopsided trees yet.

Let $c_1\leq c_2\leq\cdots\leq c_b$ be fixed positive integers. A
tree is said to be lopsided if it is $b$-ary
rooted and for each node, the edge to its $j$th child has length $c_j$, $1\leq j\leq b$.

A random lopsided tree can be constructed incrementally in the
following way:
The tree $\tilde{\mathcal{T}}_1$ consists of a single internal node \o,
the root. Additionally, we attach~$b$ external children to the root node.
If $\tilde{\mathcal{T}}_n $ already exists, take an external node
uniformly and at random and replace it by an internal node. The weights
of the edges from that internal node to its $b$ external children are
$c_1, \ldots, c_b$.
It is clear that $Z$ in the weighted $b$-ary tree framework has to be
chosen as $Z\stackrel{d}=c_W$ where $W$ is a uniform distributed
random variable on the set $\{1,\ldots, b\}$.

Then, with Corollary \ref{cor4.6.8} and Remark \ref{rem3.3.15}, by
embedding the lopsided trees in the $b$-ary tree model, and by
identifying $(T_{\tau_n})=(\tilde{\mathcal{T}}_n)$ we have the
following result for the profile $(U_l(n))$:
\begin{theorem}
If $K$ is a compact subset of
$V^*:= \{\sum_{j=1}^b c_jz^{c_j}\dvtx z \in\tilde{V}\}$ with
$\tilde V=(z_0,z_1)$ where $z_0,z_1$ are defined in Remark \ref{rem3.3.15},
then, almost surely,
\[
\lim_{n\rightarrow\infty}\sup_{l:{(b-1)l}/{\log n} \in K}
\biggl|\frac{U_{l}(n)}{\hat A_{{(b-1)l}/{\log n}}(n)}-M_{\infty}
(z_{l,n})\biggr|=0,
\]
where
$z_{l,n}\in\tilde{V}$ is the solution $z$ of $\sum
_{j=1}^bc_jz^{c_j}=\frac{(b-1)l}{\log n}$, and
\begin{eqnarray*}
\hat A_{{(b-1)l}/{\log n}}(n)
=\frac{n^{{(\sum_{j=1}^bz_{l,n}^{c_j}-1)}/{(b-1)}}}{z_{l,n}^{l}\sqrt
{2\pi\log(n){1}/{(b-1)}\sum_{j=1}^b c_j^2z_{l,n}^{c_j}}}\frac{\Gamma
({1}/{(b-1)})}{\Gamma({(\sum
_{j=1}^bz_{l,n}^{c_j})}/{(b-1)})}.
\end{eqnarray*}
We have
\[
M_{\infty}(z)=\sum_{j=1 }^b z^{c_j}\bigl(U^{(j)}\bigr)^{{1}/{(b-1)}(\sum
_{r=1}^b z^{c_r}-1)}M_{\infty,(j)}(z)
\]
with\vspace*{-2pt} $M_{\infty}(z) \stackrel{d}=M_{\infty,(j)}(z), j=1,\ldots, b$,
$U^{(j)}=\frac{Y_j}{\sum_{r=1}^b Y_r},$
 where $Y_j$ are i.i.d.\break random $\operatorname{Gamma}(\frac{1}{b-1},\frac
{1}{b-1})$ distributed random variables, and
$M_{\infty}(z)$,\break $M_{\infty,(1)}(z), \ldots, M_{\infty,(b)}(z)$,
$(U^{(1)}, \ldots,U^{(b)})$ are independent.
\end{theorem}
\end{Example}
\begin{Example}
Consider the following tree model. Start with a single internal node.
At each step the tree is expanded by choosing uniformly and at random
an internal node out of the tree and then by replacing it with a given
deterministic tree $T^*$. This model can be described by the model of
lopsided trees. Assume $|T^*|=k$.
Imagine a lopsided tree in which each replaced node gives birth to $k$
children with edge weights equal to the distances of the nodes in the
tree $T^*$ to the root of $T^*$.
The internal profile can now be calculated using the external profile
of the corresponding lopsided tree.
\end{Example}
\begin{Example}[(Plane oriented and linear recursive trees)]\label{ex3.3.5}
Plane oriented trees (PORTs) are rooted trees in which the children of
every node are oriented. The depths of nodes in random PORTs have been
studied by \citet{Mahmoud92b} and their height by \citet{Pittel94}.
PORTs can be built recursively; start with one single node, the root.
If $\tilde{\mathcal{T}}_n$ already exists, add node $v_{n+1}$ uniformly
and at random in one of the slots available. The slots are the
positions in the tree that lead to different new trees. One can think
of the slots as external nodes that are placed before, between and
after internal nodes.
So a node with $k\geq1$ children has $k+1$ external nodes attached to
it, always one external node between two (internal) children and one in
front of the first (internal) child as well as one after the last
(internal) child. If an internal node has no children, then we attach
one external node to it as a potential child.

A more general model of recursive trees is based on \citet{Pittel94}.
In these recursive trees each node $u$ has a weight $w_u$. When growing
this kind of tree, a~new node is added as a child of node $u$ with
probability proportional to $w_u$.
Now $w_u:=1+\beta\deg(u)$, where $\deg(u)$ denotes the number of
children of~$u$ and $\beta\geq0$, is called the parameter. When $\beta
$ is an integer, we can use the general tree model of \citet
{BroutinDevroye05} to describe those trees.
Let $\beta\in\mathbb{N}$ and~$\mathcal{T}^{\beta}_n$ be such a~random
recursive tree with parameter $\beta$ and with $n$ internal nodes where
$\mathcal{T}^{\beta}_1=\{\o\}$ consists of a single node, the root.
The tree is expanded by adding a child to node~$u$ with probability
proportional to $1+\beta\operatorname{deg}(u)$.
Alternatively we can choose an external node uniformly and at random
where we attached to each internal node~$u$ $\deg(u)\beta+1$ external
nodes. So when we pick an external node at level $d$ and replace it by
an internal node, we attach $\beta+2$ new external nodes to the tree,
$\beta+1$ on level $d$ and one at level $d+1$.

Now consider $(\beta+2)$-ary weighted trees $(T_t)_{t\geq0}$ where the
tree process is stopped when having $n$ internal nodes. When choosing
an external node uniformly and at random from the set of all external
nodes and when replacing it by an internal node, we add $\beta+2$
external nodes to that new internal node with weights
$Z=(Z_1,\ldots, Z_{\beta+2})$ where $Z_{(j)}=0, 1\leq j\leq\beta+1,
Z_{(\beta+2)}=1$ (the brackets in the index means that the weights are
ordered by increasing values). The external profile of that tree has a
similar distribution as the external profile of the random recursive
tree with parameter $\beta$.
Let $U_l(n)^{\beta}$ be the number of external nodes in the tree
$\mathcal{T}_n^{\beta}$ on level $l$ and $U_l(n)$ be the number of
external nodes in the tree $\mathcal{T}_n=T_{\tau_n}$, the
corresponding weighted $(\beta+2)$-ary tree. Then $U_{l+1}(n+1)^{\beta
}\stackrel{d}=U_l(n)$.
Note that for $\beta=0$ we obtain the random recursive tree of Example
\ref{ex3.3.2} and for $\beta=1$ the so called PORTs.
We can choose $Z\stackrel{d}= \operatorname{Bernoulli} (\frac{1}{\beta
+2})$ with $|\{j\leq\beta+2\dvtx Z_j=0\}|=\beta+1$.

Before formulating the convergence theorem for these recursive trees
with parameter $\beta\in\mathbb{N}$ we remark that the profile of
plane-oriented recursive trees ($\beta=1$) was analysed by \citet{Hwang05}.
For $\alpha\in[0,\frac{1}{2}]$ he obtained convergence in
distribution and of all moments of the normalized profile, where the
limit is uniquely characterized by its moment sequence.
\citet{Hwang05} presented no solution for the problem of convergence
for $\alpha\notin[0,\frac{1}{2}]$ since for $\alpha\notin[0,\frac
{1}{2}]$ only convergence of a finite number of moments is possible. As
a consequence, the characterization of the limit via moments is not
possible. In addition, no fixed point equation had been known until
now. Hwang (\citeyear{Hwang05}) anticipated convergence in distribution of
the normalized profile for $\alpha\in(\frac{1}{2}, \alpha^*)$ where
$\alpha^*$ is the solution of
$\frac{1}{2} + z + z \log(2z)=0$.

We show here even more; namely, we prove uniform almost sure
convergence for $\alpha$ in any compact subset of $(0, \alpha^*)$. Note
also that our construction also shows how the tree could be split into
subtrees in order to use the contraction method.
We identify $(\mathcal{T}^{\beta})_{n\geq0}=(T_{\tau_{n-1}}^*)_{n\geq
0}$ where the trees $(T_t^*)_{t\geq0}$ are defined as the trees
$(T_t)_{t\geq0}$ with the root node (resp. the imaginary edge $e_0$ to
the root node) having itself the weight 1.
Then we obtain
\[
\rho_t^*(l+1):=\biggl|\biggl\{u \in\partial T_t^*\dvtx D_u=\sum_{e\in\pi(u)}Z_e=l+1\biggr\}
\biggr| =\rho_t(l).
\]
Note that $e_0\in\pi(u)$ for all $u$.
Finally, it follows that
\[
U_{l+1}(n+1)^{\beta}=\rho_{\tau_{n}}^*(l+1)=\rho_{\tau_n}(l).
\]
\begin{theorem}
Let $K \in(0,z_0)$ be a compact set where $z_0$ is the only solution
of $z\log(z)-z-\beta=0$. Then almost surely
\[
\lim_{n\rightarrow\infty}\sup_{l:{(\beta+1)l}/{\log n }\in K}
\biggl|\frac{U_{l_n+1}(n+1)^{\beta}}{\hat A_{{(b-1)l}/{\log
n}}(n)}-M_{\infty}(z_{l,n})\biggr|=0
\]
with
\[
\hat A_{{(\beta+1)l}/{\log n }}(n)=\frac{n^{{(\beta
+z_{l,n})}/{(\beta+1)}}}{z_{l,n}^{l}\sqrt{2\pi l}}\frac{\Gamma(
{1}/{(\beta+1)})}{\Gamma(1+{z_{l,n}}/{(\beta+1)})},
\]
where $z_{l,n}:=(\beta+1)\frac{l}{\log n}$.
Further
\[
M_{\infty}(z)=\sum_{j=1}^{\beta+1}\bigl(U^{(j)}\bigr)^{{(\beta
+z)}/{(\beta+1)}}M_{\infty, (j)}(z)+\Biggl(1-\sum_{j=1}^{\beta
+1}U^{(j)}\Biggr)z M_{\infty, (\beta+2)}(z),
\]
where
$M_{\infty, (1)}(z), \ldots,M_{\infty, (\beta+2)}(z) $ are
independent, distributed as $M_{\infty}(z)$\break and independent of
$(U^{(1)},\ldots, U^{(\beta+1)})$ with $U^{(j)}=\frac{Y_j}{\sum
_{r=1}^{\beta+2}Y_r}$\vspace*{-3pt} where\break $Y_j\stackrel{d}=\operatorname{Gamma}
(\frac{1}{\beta+1},\frac{1}{\beta+1})$ i.i.d.
\end{theorem}
\begin{Remark}
For $\beta=1$, $z_0$ is the only solution of $\frac{1}{2} + \frac
{z}{2}-\frac{z}{2}\log(z)=0$. So~$l$ has to be chosen so that $\frac
{l}{\log n} \in(0,z_0/2)$. Obviously $z_0/2$ is the only solution of
$\frac{1}{2} + z-z\log(2z)=0$, $z_0/2=\alpha^*$.
\end{Remark}

\begin{pf}
We have $b=\beta+2$, $Z\stackrel{d}= \operatorname{Bernoulli} (\frac{1}{\beta
+2})$, $bEz^Z=\beta+1+z$.
It follows that
\begin{eqnarray*}
\tilde{V}
& = & \{z \geq0\dvtx  1-bEz^Z<-\log(z)bE(Zz^Z)\}\\
& = & \{z \geq0\dvtx  z \log(z)-z-\beta<0\}.
\end{eqnarray*}
Since $bE(Zz^Z)=z$ we have $\tilde{V}=V^*$.
Define $f(z):= z \log(z)-z-\beta$. Then
\[
f'(z)=\log(z), \qquad  f''(z)>0,
\]
and $(1,-1-\beta)$ is a local minimum of $f$.
Since
\[
\lim_{z \searrow0}f(z)=-\beta, \qquad \lim_{z\rightarrow\infty
}f(z)=\infty,
\]
there is only one solution of $f(z)=0$ that we will call $z_0$.
In the interval $[0,z_0)$ $f$ is negative other than that nonnegative.
The rest follows from Corollary \ref{cor4.6.8}.
\end{pf}

 Note that independently of this work \citet{Sulzbach08} proved a
functional limit theorem for the profile of plane oriented recursive
trees using the martingale method.
\end{Example}
\begin{Example}[(Changes of direction in a binary search tree)]
Let $\tilde{\mathcal{T}}_n$ be a random binary search tree with $n$
internal nodes, and let $u \in\tilde{\mathcal{T}}_n$.
Define $D_n(u):=D_n(\pi(u))$ as the number of changes of direction in
$\pi(u)$ where $\pi(u)$ is the path from the root to node $u$.
Now let $0$ and $1$ encode a move down to the left and to the right,
respectively.
For example the path encoded by $1001010110$ will have $D=7$, that is,
a count of each occurrence of the patterns $01$ and $10$.

We are interested in $D_l(n):=|\{u \in \partial\tilde{ \mathcal
{T}}_n:  D_n(u)=l\}|.$
Broutin and Devroye~(\citeyear{BroutinDevroye05}) introduced the following
labelling of the edges: for each level $l\geq2$ of edges form the word
$(0110)^{l-1}$, and map the binary characters to the edges from left to
right. Call this weighted binary tree $\mathcal{T}_n$. Then, by
embedding, we find that $D_l(n)=|\{u \in\partial{\mathcal{T}}_n\dvtx
D_u=l\}|$.
Consequently choose $Z\stackrel{d}=\operatorname{binomial}(\frac{1}{2}
)$, $Z_2=1-Z_1$, and obtain the following as in the random recursive
tree case:
\begin{theorem}
Let $K\subset(0,e)$ be a compact subset. Then almost surely
\[
\lim_{n\rightarrow\infty}\sup_{l:{l}/{\log n}\in K}\biggl|\frac
{D_{l}(n)}{ED_{l}(n)}-M_{\infty}\biggl(\frac{l}{\log
n}\biggr)\biggr|=0,
\]
where
\[
M_{\infty}(c)\stackrel{d}=c
U^{c}M_{\infty}(c)+(1-U)^{c}M_{\infty}^*(c),
\]
where $M_{\infty}(c)\stackrel{d}=M_{\infty}^*(c)$, $U$ is a
$\operatorname{uniform }[0,1]$ random variable and
$M_{\infty}(c)$, $M_{\infty}^*(c)$ and $U$ are independent.
\end{theorem}
\end{Example}
\begin{Example}[(Random $l$-colouring of the edges in a tree)]
Take a random binary search tree and randomly color the edges with one
of~$l$ different colors.
We can think of different problems in that framework. For instance, we
could be interested in the question how many nodes $u$ have exactly~$l$
red edges in $\pi(u)$ if color red appears with probability $p$.
For this problem we have to choose $Z\stackrel{d}= \operatorname{Bernoulli} (p)$,
since we count only red edges which we mark with $Z=1$ and all other
colored edges with $Z=0$.
Let $D_n(l)$ be the number of nodes in the tree $\mathcal{T}_n$ with
exactly $l$ red edges in $\pi(u)$.
\begin{theorem}
Let $K\subset\{2zp\dvtx z\in\tilde{V}\}$ with $\tilde{V}=\{z>0\dvtx 2pz\log
(z)-2zp+2p-1<0\}$.
Then almost surely
\[
\lim_{n\rightarrow\infty}\sup_{l:{l}/{\log n}\in K}\biggl|\frac
{D_n(l)}{\hat A_{{l}/{\log n}}(n)}-M_{\infty}(z_{l,n}
)\biggr|=0,
\]
where we set $z=z_{l,n}:=\frac{l}{2p\log n}$ and
\begin{eqnarray*}
\hat A_{{l}/{\log n}}(n)&:=&\frac{n^{1-2p+2pz_{l,n}}}{z_{l,n}^{l}\sqrt
{2\pi l} }\frac{1}{\Gamma(2(1-p)+2pz_{l,n}))},
\\
M_{\infty}(z)&\phantom{:}=&z^{Z_1}U^{1-2p+2pz}M_{\infty,
(1)}(z)+z^{Z_2}(1-U)^{1-2p+2pz}M_{\infty, (2)}(z),
\end{eqnarray*}
where $M_{\infty, (1)}(z),M_{\infty, (2)}(z) $ are independent,
distributed as $M_{\infty}(z)$, independent of $Z_1,Z_2,U$ where
$U\stackrel{d}=\operatorname{uniform}[0,1]$ and $Z_1,Z_2$ are independent,
identically distributed with  $\operatorname{Bernoulli} (p)$ distribution.
\end{theorem}
\begin{pf}
Since $b=2,  bEZz^Z=2zp, bEz^Z=2(1-p+zp)$, we have
$\tilde{V}=\{z>0\dvtx 2pz\log(z)-2zp+2p-1<0\}$ and $V^*=\{2zp\dvtx z\in\tilde
{V}\}$.
\end{pf}

For the random recursive tree the number of nodes with paths having
exactly~$l$ red edges can be analyzed taking $Z\stackrel{d}=
\operatorname{Bernoulli} (p)\times \operatorname{Bernoulli} (1/2) $, thus having $Z\stackrel
{d}= \operatorname{Bernoulli} (p/2)$ [cf. Example \ref{ex3.3.2}]. The random
recursive tree can be interpreted as a weighted binary tree. Now
randomly color the edges of this tree. The probability of having a red
edge is now $p$ and the probability of having an edge with weight $1$
is $1/2$.
This model could alternatively be analyzed in a~2-dimensional weighted
model where
$Z=(Z^{(1)},Z^{(2)})$, $Z^{(1)}\stackrel{d}= \operatorname{Bernoulli} (1/2)$,
$Z^{(2)}\stackrel{d}= \operatorname{Bernoulli} (p)$
where $Z^{(1)},Z^{(2)}$ are independent.
\end{Example}
\begin{Example}[(The left minus right exceedance)]
Let $\mathcal{T}_n$ be a binary search tree with $n$ internal nodes and
let $u \in\mathcal{T}_n$.
Define $D_u:=\sum_{e \in\pi(u)}(L(e)-R(e))$ where $L(e)$ is the
indicator function of $e$ being a left edge and $R(e)$ is analogously
the indicator of $e$ being a right edge.
We are interested in
\[
U_l(n):=|\{u \in\partial{\mathcal{T}_n}\dvtx  D_u=l\}|,
\]
namely, the number of external nodes in our binary search tree which
have exactly~$l$ more left edges than right edges in the path from the
root leading to that external node. Naturally in the framework of
weighted $b$-ary trees we choose $b=2$ and mark all right edges with
$-1$ and weight all left edges with $1$. Since right and left edges are
equally likely to be chosen, we use $Z$ with $P(Z=1)=\frac{1}{2}=P(Z=-1)$.
%
\begin{theorem}
Let $K$ be a compact subset of $ (z_0-\frac{1}{z_0},z_1-\frac{1}{z_1})$
where $0<z_0<z_1$ are the two positive solutions of
\[
1-z-\frac{1}{z}+z\log(z)-\frac{1}{z}\log(z)=0.
\]
Then almost surely,
\[
\lim_{n\rightarrow\infty}\sup_{l: {l}/{\log n} \in K}\biggl(\frac
{U_{l}(n)}{\hat A_{{l}/{\log n}}(n)}-M_{\infty}(z_{l,n})\biggr)=0,
\]
where we define
\[
\hat A_{l/\log n}(n)
:=\frac{1}{\Gamma(z_{l,n}+{1}/{z_{l,n}})}\frac
{n^{z_{l,n}+{1}/{z_{l,n}}-1}}{z_{l,n}^{l}\sqrt{2\pi(z_{l,n}+
{1}/{z_{l,n}})\log(n)}}
\]
with
\[
z_{l,n}:=\frac{l}{2\log n}+\sqrt{\biggl(\frac{l}{2\log n}\biggr)^2+1}.
\]
Further we have
\[
M_{\infty}(z)= z(U)^{z + {1}/{z}-1}M_{\infty, (1)}(z)
+ \frac{1}{z}(1-U)^{z + {1}/{z}-1}M_{\infty, (2)}(z),
\]
where
$M_{\infty, (1)}(z),M_{\infty, (2)}(z) $ are independent, distributed
as $M_{\infty}(z)$ and independent of $U\stackrel{d}= \operatorname{uniform}[0,1]$.
\end{theorem}
\begin{pf}
First note that
\begin{eqnarray*}
\tilde{V}&=&\{z>0\dvtx 1-bEz^Z+\log(z)bE(Zz^Z)<0\}{}\\
&=& \biggl\{z>0\dvtx 1-\biggl(z+\frac{1}{z}\biggr)+\log(z)\biggl(z-\frac
{1}{z}\biggr)<0\biggr\}.
\end{eqnarray*}
Easily we obtain for $z>0$ with $g(z):=1-(z+\frac{1}{z})+\log(z)(z-\frac{1}{z})$
\begin{eqnarray*}
g'(z)&=&\log(z)\biggl(1+\frac{1}{z^2}\biggr)=0\quad \Leftrightarrow\quad  z=1, \\
g''(1)&>&0,
\end{eqnarray*}
that $g$ has a single local minimum at $z=1$ with $g(1)=-1$
and since $g$ is continuous on $(0,\infty)$ with
\[
\lim_{z\searrow0}g(z)=\infty, \qquad  \lim_{z \rightarrow\infty
}g(z)=\infty,
\]
there are exactly two roots of the function $g$ on $(0,\infty)$.
Call them $0<z_0<z_1$.
Now $V^*:=\{z-\frac{1}{z}\dvtx  z \in(z_0, z_1)\}=(z_0-\frac{1}{z_0},
z_1-\frac{1}{z_1})$.
If $c:=\frac{l}{\log(n)}$, then choose $z=z(c)>0\dvtx  z-\frac{1}{z}=c
\Leftrightarrow z^2-cz-1=0 \Leftrightarrow z=\frac{c}{2}+\sqrt{
(\frac{c}{2})^2+1}.$
From this the proof follows.
\end{pf}
\end{Example}
\begin{Example}[(Stochastic models for the web graph)]
We give a new example not contained in \citet{BroutinDevroye05}.
The web may be viewed as a~directed graph in which each vertex is a
static HTML web page, and each edge is a hyperlink from one web page to
another. \citet{Kumar00} proposed and analyzed a class of random graph
models inspired by empirical observations on the web graph. These
observations suggested that the web is not well modeled by traditional
graph models.

The linear growth copy model of \citet{Kumar00} is parameterized by a
copy factor $\alpha\in(0,1)$ and a constant outdegree $d \geq1$.
Only the choice $d=1$ results in a random forest that might be turned
into a tree and studied using our framework.

We start with one single vertex.
Assume that the random forest $\tilde{\mathcal{T}}_n$ with $n$ internal
nodes has already been created.
At each time step, one vertex $u$ is added by the following procedure:
from the tree $\tilde{\mathcal{T}}_n$ choose a vertex uniformly and at
random. Call this vertex $v$. With probability $\alpha$ we attach node
$u$ as a child to node $v$. With probability $1-\alpha$ the node $u$
becomes a brother of node $v$; that means that we attach node $u$ as a
child to the father of node $v$. If node $v$ is a root with no
ancestors, we let $u$ be an isolated node, namely the root of a new
tree consisting of that single node.
We could now ask how many nodes are roots, nodes on level $1,2,\ldots$
and so on in that random forest.

We can interpret this random forest as a binary tree with weighted
edges. When raising the forest we may instead raise the binary tree as follows.
In the random forest a new node $u$ is attached by choosing uniformly
and at random an internal node $v$ out of the existing forest $\tilde
{\mathcal{T}}_n$. In the binary tree $\mathcal{T}_{n-1}=T_{\tau_{n-1}}$
we will instead choose an external node, call it $\tilde v$, uniformly
and at random from one of the~$n$ external nodes. With the probability
$\alpha$, the new node in the forest will be a child of node~$v$ and
located one level below $v$. We transmit this by making the external
node~$\tilde{v}$ in the binary tree an internal one and attach two new
external nodes to $\tilde v$, one with edge weight $0$, representing
$v$, and the other with edge-weight $1$, representing $u$.
With the probability $1-\alpha$, node $u$ becomes a brother of $v$,
that means it stays on the same level as node $v$.
Then we will replace the external node $\tilde{v}$ in the binary tree
by an internal node and attach two new external nodes to it, one with
edge weight $0$, representing $v$, and the other with edge-weight $0$,
representing~$u$.
Then an arbitrary edge has weighted one with probability $\alpha/2$ and
otherwise it has weight zero.\vspace*{-1pt}
Choose $Z\stackrel{d}= \operatorname{Bernoulli} (\frac{\alpha}{2})$
and the weights $Z_1,Z_2$ attached to the root of the binary tree as follows.\vspace*{-1pt}
Let $Z_1\stackrel{d}=Z$ and $Z_2=\mathbh{1}_{\{Z_1=0\}}Y$ with
$Y\stackrel{d}= \operatorname{Bernoulli} (\frac{\alpha}{2-\alpha})$
and being independent of $Z_1$. Then $Z_2\stackrel{d}=Z$ and the
resulting tree $\mathcal{T}_{n-1}$ grows as the tree described above.
By embedding we obtain the following:
\begin{theorem}
Let $K$ be a compact subset of $V^*:=(\alpha z_0,\alpha z_1) $ with
$z_0, z_1$ the two roots of the function $f(z)=-1+\alpha+\alpha
z(1+\log(z))$.
Then almost surely
\[
\lim_{n\rightarrow\infty}\sup_{l: {l}/{\log n}\in K}\biggl|\frac
{U_{l}(n)}{\hat A_{{l}/{\log n}}(n)}-M_{\infty}\biggl(\frac
{l}{\alpha\log n}\biggr)\biggr|=0
\]
with
\[
\hat A_{{l}/{\log n}}(n):=\frac{n^{1-\alpha+{l}/{\log
n}}}{({l}/{(\alpha\log n)})^{l}\sqrt{2\pi l}}\frac
{1}{\Gamma(2-\alpha+{l}/{\log n})}
\]
and
\[
M_{\infty}(z)=z^{Z_1}(U)^{1-\alpha+c}M_{\infty, (1)}(z)+
z^{\mathbh{1}_{\{Z_1\}}Y}(1-U)^ {1-\alpha+c} M_{\infty, (2)}(z),
\]
where $M_{\infty, (1)}( z), M_{\infty, (2)}( z) $ are independent,
distributed as $M_{\infty}(z)$ and independent of $U$ where $U\stackrel
{d}=\operatorname{uniform}[0,1]$.
\end{theorem}
\end{Example}

\begin{Example}[(Combination of weights)]
Higher dimensional weights can be used to describe all the trees
studied earlier with additional weights attached to the nodes or,
alternatively, edges.
For example we can study a 2-ary tree with $Z=(Z^{(1)},Z^{(2)})$,
$Z^{(1)}=1$ and $Z^{(2)}\stackrel{d}=$ binomial$(\frac{1}{2})$ which\vspace*{1pt}
refers to a random binary search
tree with edges marked with zero or one.
We can think of a situation where we use the second weight for
identifying if the ancestor passes some attribute on to its child
($=1$) or not ($=0$). Let
\[
U_{(\tilde l,l)}(n)=\biggl|\biggl\{u \in\partial{ \mathcal{T}}_n\dvtx
D_u=\biggl(\sum_{e \in\pi(u)}Z^{(1)}_e, \sum_{e \in\pi
(u)}Z^{(2)}_e\biggr)=(\tilde l,l)\biggr\}\biggr|.
\]
Applying Theorem \ref{theo2.3.9} (resp. Corollary \ref{cor3.3.13}) we
can then describe the asymptotics of that numbers.
\end{Example}
%


\printaddresses

\end{document}